\documentclass[preprint,12pt]{elsarticle} 
\biboptions{numbers,sort&compress}
\usepackage[colorlinks=true, linkcolor=blue, citecolor=blue, urlcolor=blue]{hyperref}
\usepackage{algorithm}
\usepackage{algpseudocode} 
\usepackage{amsmath, amssymb}
\usepackage{graphicx}
\usepackage{subfig}
\usepackage{float}
\usepackage{tikz}
\usepackage{pgfplots}
\usepackage{caption}
\usepackage{url}
\usepackage{booktabs}

\journal{Journal of Computational Physics}

\begin{document}

\begin{frontmatter}

\title{A finite element method using a bounded auxiliary variable for solving the Richards equation}

\author[aff1]{Abderrahmane Benfanich}
\ead{abenf099@uottawa.ca}

\author[aff1]{Yves Bourgault}
\ead{ybourg@uottawa.ca}

\author[aff1,aff2]{Abdelaziz Beljadid\corref{cor1}}
\ead{Abdelaziz.BELJADID@um6p.ma}
\ead{abeljadi@uottawa.ca}

\address[aff1]{Department of Mathematics and Statistics, University of Ottawa, Canada}
\address[aff2]{Mohammed VI Polytechnic University, Morocco}

\cortext[cor1]{Corresponding author.}

\begin{abstract}
The Richards equation, a nonlinear elliptic parabolic equation, is widely used to model infiltration in porous media. We develop a finite element method for solving the Richards equation by introducing a new bounded auxiliary variable to eliminate unbounded terms in the weak formulation of the method. This formulation is discretized using a semi-implicit scheme and the resulting nonlinear system is solved using Newton's method. Our approach eliminates the need of regularization techniques and offers advantages in handling both dry and fully saturated zones. In the proposed techniques, a non-overlapping Schwarz domain decomposition method is used for modeling infiltration in layered soils. We apply the proposed method to solve the Richards equation using the Havercamp and  van Genuchten models for the capillary pressure. Numerical experiments are performed to validate the proposed approach, including tests such as modeling flows in fibrous sheets where the initial medium is totally dry, two cases with fully saturated and dry regions, and an infiltration problem in layered soils. The numerical results demonstrate the stability and accuracy of the proposed numerical method. The numerical solutions remain positive in the presence of totally dry zones. The numerical investigations clearly demonstrated the capability of the proposed method to effectively predict the dynamics of flows in unsaturated soils.
\end{abstract}


\end{frontmatter}

\section{Introduction}
\label{sec:introduction}
Studying fluid dynamics within porous media is critical to numerous scientific and engineering disciplines, encompassing hydrogeology, environmental engineering and soil science. These fields rely heavily on understanding infiltration to address such as ranging from sustainable water resource management and contamination mitigation. 


Our interest specifically lies in the infiltration process, which describes the movement of water through unsaturated porous media. A commonly used model for infiltration through unsaturated soils is the Richards equation, which was obtained from a multiphase flow extension of Darcy's law and the principle of mass conservation \cite{Richards1931}. The versatility of the Richards equation extends its applicability beyond classical soil hydrology to material science and engineering. However, accurately modeling domains that include both completely dry ($S=0$) and fully saturated ($S=1$) regions remains challenging, as existing formulations have difficulty addressing both cases simultaneously. This highlights the necessity of developing robust formulations capable of seamlessly modeling these phenomena without requiring variable switching or encountering numerical instabilities.

The Richards equation has multiple formulations: the pressure-based formulation ($\Psi$), the saturation-based formulation ($S$) and the mixed formulation ($S, \Psi$), as detailed in \cite{Celia1990}. Each formulation has its advantages and drawbacks. The $\Psi$-based formulation can be applied to both fully saturated and partially unsaturated flows (where saturation approaches but does not reach zero). As shown in \cite{AltLuckhaus1983}, the pressure head formulation has a unique and continuous solution. The $\Psi$-based form has been frequently used in various applications, for instance, in \cite{Haverkamp1977}. However, it has been demonstrated in \cite{Celia1990} that the $\Psi$-based formulation may produce significant global mass balance errors unless very small time steps are employed. Some of these difficulties are mitigated by mixed-form schemes, which provide improved mass conservation properties. In \cite{Celia1990}, the mixed formulation was addressed by employing a modified Picard iteration scheme.

Saturation-based formulations exhibit better numerical performance compared to the $\Psi$-based methods, especially in very dry soils, as observed in \cite{Hills1989}. Additionally, it has been proven in \cite{Marinoschi2006} that the saturation-based form has a unique solution when the medium is not fully saturated. However, these formulations suffer from substantial numerical difficulties. These challenges primarily stem from the discontinuous nature of saturation as a primary variable and from complications in saturated regions, where soil-water diffusivity may become unbounded in certain capillary and saturation constitutive relationships.

A method used to overcome these difficulties is the primary variable switching approach, in which either pressure head or saturation is selected as the primary variable, depending on local saturation conditions at numerical nodes, as seen in \cite{Diersch1999, Forsyth1995, Wu2001}. Typically, this method is implemented within Newton-Raphson iterative schemes, where the Jacobian matrix is formed based on derivatives with respect to the chosen primary variable (e.g., \cite{Brunner2012,Krabbenhoft2007}). Despite its general applicability, variable-switching approaches often introduce non-smooth transitions between primary variables, which can lead to physically unrealistic solutions, especially at the interface between saturated and unsaturated zones, a problem noted in \cite{Krabbenhoft2007, Zha2017}. Techniques to address these numerical challenges, such as smoothing techniques and improved switching criteria, have been developed, but their effectiveness remains problem-specific, as discussed in \cite{Maina2017, Kees2002, Krabbenhoft2007}. Additionally, these techniques are primarily applicable within Newton-Raphson iterative schemes. A review of these methods is provided in \cite{Zha2019}.

Another strategy to deal with the degeneracy and nonsmoothness of the Richards equation is the use of regularization techniques. In this approach, degenerate constitutive relationships are replaced by smoothed counterparts that remain well-posed, controlled by a regularization parameter. Various schemes have been proposed in the literature, ranging from parabolic regularizations for unsaturated flow with dry regions \cite{Schweizer2007} to adaptive regularization frameworks guided by a posteriori error estimators \cite{Fevotte2024}, as well as schemes tailored for doubly degenerate equations and outflow conditions \cite{Pop2011}. These methods improve robustness and ensure convergence of numerical solvers, though they may introduce an additional modeling parameter that should be selected.

To address these limitations, we introduce a new formulation based on a single primary variable $u$, eliminating the need for variable switching. Unlike saturation-based formulations, our approach remains well-defined in both fully unsaturated regions ($S=0$) and fully saturated regions ($S=1$). However, similar to the $\Psi$-based formulation, our method initially exhibits difficulties related to the conservation of global mass. We overcome this challenge by employing a semi-implicit iteration to solve the mixed formulation in the variables ($u, S$). The resulting nonlinear system is then solved using Newton’s method, where $u$ is the main variable.

To effectively handle heterogeneity, domain decomposition techniques have been employed, as described in \cite{Lion1990}. In this context, the Richards equation is solved separately within each distinct soil type, introducing appropriate transmission conditions at the interfaces between subdomains. These conditions exploit the continuity of both water flux and the pressure head ($\Psi$) across interfaces. This approach has proven effective in accurately solving scenarios involving heterogeneous soils with combined saturated and unsaturated regions. Its performance has been demonstrated by modifying the well-known Manzini test case from \cite{Manzini2004} to include saturated-unsaturated heterogeneous layered soils.


We propose a numerical method based on the \(u\)-formulation to solve the Richards equation. 
This approach employs a finite element discretization in space, in line with previous works on Galerkin finite element methods \cite{Arbogast1993, Slodicka2002, Bergamaschi1999}, and incorporates a time-stepping strategy related to semi-implicit schemes \cite{Keita2021} where we solve the resulting nonlinear system by Newton's method. 

The article is organized as follows. We begin in Section \ref{sec:waterinfiltration} by introducing the Richards equation, detailing its various formulations and reviewing empirical models for hydraulic functions. This comprehensive section also motivates the derivation of a new $u$-formulation, designed to overcome the limitations of existing approaches and presents the foundational equations for the domain decomposition method used to effectively address heterogeneous media problems. Following this, we delve into the variational formulations for the pressure-based ($\Psi$), saturation-based ($S$) and our novel $u$-formulations, enabling a clear comparison of their mathematical structures. Subsequently, in Section \ref{sec:NMETHODS}, we outline the numerical methods employed for solving these equations, including temporal and spatial discretization schemes, along with specialized approaches for the $u$-formulation and the domain decomposition method. Section \ref{sec:Numres}, presents numerical results that validate our proposed method across various challenging scenarios, from fully unsaturated fibrous sheets to heterogeneous test cases with both saturated and dry regions, where we calculate the order of convergence using a manufactured solution approach. Finally, in Section \ref{sec:conclusion}, we summarize the main findings of this work and highlight directions for future research, including rigorous numerical analysis of the proposed method and extensions to coupled processes such as solute and heat transport.

\section{Notation and Functional Spaces}

In this paper, we use the following notation and definitions for functional spaces and mathematical symbols:

\begin{itemize}
\item $\Omega \subset \mathbb{R}^d$, with $d=1,2,3$, denotes a bounded spatial domain with boundary $\partial \Omega$.

\item $I = (0,T)$ represents the temporal domain with final time $T > 0$.

\item $L^2(\Omega)$ is the space of square-integrable functions on the domain $\Omega$, equipped with the norm

$$
\|u\|_{L^2(\Omega)} = \left(\int_{\Omega}|u|^2\,dx\right)^{1/2}.
$$

\item $H^1(\Omega)$ denotes the standard Sobolev space of functions with square-integrable first derivatives, defined as

$$
H^1(\Omega) = \{u \in L^2(\Omega): \nabla u \in (L^2(\Omega))^d\},
$$

with norm

$$
\|u\|_{H^1(\Omega)} = \left(\|u\|_{L^2(\Omega)}^2 + \|\nabla u\|_{L^2(\Omega)}^2\right)^{1/2}.
$$

\item $H_{\Gamma_D}^1(\Omega)$ denotes the subspace of $H^1(\Omega)$ satisfying homogeneous Dirichlet boundary conditions on $\Gamma_D \subset \partial \Omega$:

$$
H_{\Gamma_D}^1(\Omega) = \{u \in H^1(\Omega): u|_{\Gamma_D}=0\}.
$$

\item $L^2(I;X)$ is the Bochner space of functions defined on the interval $I$ taking values in a Banach space $X$, with norm

$$
\|u\|_{L^2(I;X)} = \left(\int_I \|u(t)\|_X^2\,dt\right)^{1/2}.
$$

\item The variable $\boldsymbol{x}$ denotes the spatial coordinate vector.

\item $\frac{\partial}{\partial \boldsymbol{x}} := \left( \frac{\partial}{\partial x_1}, \frac{\partial}{\partial x_2}, \dots, \frac{\partial}{\partial x_d} \right)$ denotes the partial gradient with respect to $\boldsymbol{x}$.

\item $\boldsymbol{e}_z$ denotes the unit vector in the vertical direction (positive upward).

\end{itemize}

\section{The mathematical model of infiltration}
\label{sec:waterinfiltration}

This section describes the Richards equation, the governing model for water infiltration in unsaturated porous media. We present several formulations of this equation and review empirical models for the hydraulic functions. The formulations are discussed with emphasis on their respective strengths and limitations, motivating a new formulation whose derivation is detailed. To effectively address layered soil problems, the equations for the domain decomposition method are also presented.

\subsection{The Richards equation}

Let $\Omega \subset \mathbb{R}^d$, with $d = 1, 2, 3$, be a bounded open set and let $T > 0$. Define $I = (0, T)$. The \textit{Richards equation} \citep{Richards1931} is given by
\begin{equation}
\label{eq:richards}
\frac{\partial \theta}{\partial t} + \nabla \cdot \boldsymbol{q} = 0, \quad \text{in } \Omega \times I,
\end{equation}
where $\theta = \theta(\boldsymbol{x}, t)$ is the water content and the water flux $\boldsymbol{q}$ is described by Darcy-Buckingham law \citep{Richards1931}:
\begin{equation}
\label{eq:darcy}
\boldsymbol{q} = -K_s(\boldsymbol{x}) K_r(\boldsymbol{x}, S) \nabla (\Psi(\boldsymbol{x}, S) + z).
\end{equation}
Here, $S = (\theta - \theta_r)/(\theta_s - \theta_r)$ is the effective saturation, $K_s$ is the saturated hydraulic conductivity, $K_r$ is the relative permeability and $\boldsymbol{x} = (x, z)$ denotes spatial coordinates. The vertical coordinate $z$ is oriented positively upward. The functions $\theta_s(\boldsymbol{x})$ and $\theta_r(\boldsymbol{x})$ represent the saturated and residual water content, respectively, with
\begin{equation*}
\theta_r(\boldsymbol{x}) \leq \theta(\boldsymbol{x}, t) \leq \theta_s(\boldsymbol{x}), \quad \forall (\boldsymbol{x}, t) \in \Omega \times I.
\end{equation*}

The Richards equation admits multiple formulations:
\begin{itemize}
    \item \textbf{Pressure-based formulation:} Find \( \Psi \) such that
    \begin{equation}
    \label{pressurebased}
    \phi(\boldsymbol{x}) \frac{\partial S}{\partial \Psi}(\boldsymbol{x}, \Psi) \frac{\partial \Psi}{\partial t}
    - \nabla \cdot \left( K_s(\boldsymbol{x}) K_r(\boldsymbol{x}, \Psi) \nabla (\Psi + z) \right) = 0.
    \end{equation}
    \item \textbf{Saturation-based formulation:} Find \( S \) such that
    \begin{equation}
    \label{saturationbased}
    \phi(\boldsymbol{x})\,\frac{\partial S}{\partial t}
            \;-\;
            \nabla\!\cdot\!\Bigl[
            K_{s}(\boldsymbol{x})\,K_{r}(\boldsymbol{x},S)\,
            \bigl(
            \tfrac{\partial\Psi}{\partial S}(\boldsymbol{x},S)\,\nabla S
            \;+\;
            \tfrac{\partial\Psi}{\partial\boldsymbol{x}}(\boldsymbol{x},S)
            \;+\;
            \nabla z
            \bigr)
            \Bigr]
            =
            0.
    \end{equation}
    \item \textbf{Mixed formulation :} Find \( \Psi \) and \( S \) such that
    \begin{equation*}
    \phi(\boldsymbol{x}) \frac{\partial S}{\partial t}
    - \nabla \cdot \left( K_s(\boldsymbol{x}) K_r(\boldsymbol{x}, \Psi) \nabla (\Psi + z) \right) = 0.
    \end{equation*}
\end{itemize}
Here, $\phi = \theta_s - \theta_r$ and the relation between $S$ and $\Psi$ is given by \eqref{eq:leverett}.

To close the system, the relationships between $K_r$ and $\Psi$ as functions of $S$ must be specified. Several empirical models exist in the literature, including those by \cite{Gardner1958}, \cite{Brooks1966}, \cite{Haverkamp1977} and \cite{VanGenuchten1980}.

The pressure head is related to effective saturation by Leverett’s function \citep{Leverett1941}:
\begin{equation}
\Psi(S) = h_{\text{cap}}(\boldsymbol{x}) J(\boldsymbol{x}, S),
\label{eq:leverett}
\end{equation}
where $h_{\text{cap}}: \Omega \to \mathbb{R}^+$ is the capillary rise function and $J: \Omega \times (0, 1) \to (-\infty, 0]$ is strictly monotone and $C^1$ in its second argument.

\begin{table}[H]
  \centering
  \begin{minipage}{0.9\textwidth}
    \centering
    \caption{Soil hydraulic functions for different models.}
    \label{tab:Leverett function}
    \begin{tabular}{|c|c|c|c|}
      \hline
      Model & $S(\Psi)$ & $K_r(S)$ & $J(S)$ \\
      \hline
      Gardner (1958)\footnotemark[1]
        & $e^{\alpha \Psi}$ & $S$ & $\ln(S)$ \\
      \hline
      Brooks-Corey (1966)\footnotemark[1]\footnotemark[2]
        & $(\alpha \Psi)^{-\lambda}$
        & $S^{B}$
        & $-S^{-\frac{1}{\lambda}}$ \\
      \hline
      Haverkamp (1977)\footnotemark[1]
        & $\frac{1}{1+|\alpha \Psi|^\beta}$
        & $\frac{1}{1+ |A \Psi(S)|^\gamma}$
        & $-(S^{-1}-1)^{1/\beta}$ \\
      \hline
      van\,Genuchten (1980)\footnotemark[1]\footnotemark[3]
        & $\frac{1}{(1+|\alpha \Psi|^n)^m}$
        & $S^{0.5}(1-(1-S^{1/m})^m)^2$
        & $-(S^{-1/m}-1)^{1/n}$ \\
      \hline
    \end{tabular}
    \footnotetext[1]{The parameters $\alpha,\lambda,\beta,\gamma,A,B,n$ are empirical and depend on type of soils. They satisfy: $\alpha>0,\;\lambda>0,\;\beta>1,\;\gamma>0,\;A>0,\; B>0,\;n>1,\;m=1-1/n$.}
    \footnotetext[2]{We use a power law for the relative permeability}
    \footnotetext[3]{The van Genuchten relative permeability is based on the Mualem model \citep{mualem1976}.}
  \end{minipage}
\end{table}

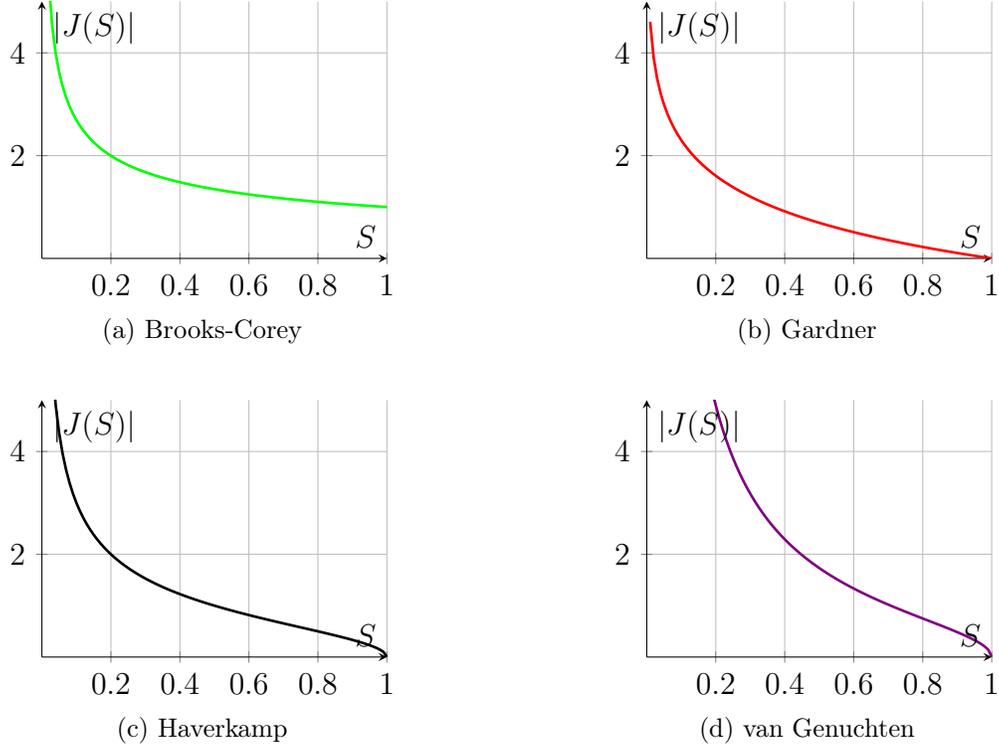
\begin{figure}[H]
\centering
\subfloat[Brooks-Corey]{
\begin{tikzpicture}
\begin{axis}[
    axis lines=middle,
    xlabel=$S$, ylabel=$|J(S)|$,
    domain=0.01:1, samples=100,
    grid=major, width=0.45\textwidth, height=5cm,
    ymin=0, ymax=5, xmin=0, xmax=1
]
\addplot[line width=1, color=green]{(x^(-0.43))};
\end{axis}
\end{tikzpicture}
}
\hfill
\subfloat[Gardner]{
\begin{tikzpicture}
\begin{axis}[
    axis lines=middle,
    xlabel=$S$, ylabel=$|J(S)|$,
    domain=0.01:1, samples=100,
    grid=major, width=0.45\textwidth, height=5cm,
    ymin=0, ymax=5, xmin=0, xmax=1
]
\addplot[line width=1, color=red]{abs(ln(x))};
\end{axis}
\end{tikzpicture}
}
\\
\vspace{0.3cm}
\subfloat[Haverkamp]{
\begin{tikzpicture}
\begin{axis}[
    axis lines=middle,
    xlabel=$S$, ylabel=$|J(S)|$,
    domain=0.01:1, samples=100,
    grid=major, width=0.45\textwidth, height=5cm,
    ymin=0, ymax=5, xmin=0, xmax=1
]
\addplot[line width=1, color=black]{(x^(-1)-1)^(1/2)};
\end{axis}
\end{tikzpicture}
}
\hfill
\subfloat[van Genuchten]{
\begin{tikzpicture}
\begin{axis}[
    axis lines=middle,
    xlabel=$S$, ylabel=$|J(S)|$,
    domain=0.01:1, samples=100,
    grid=major, width=0.45\textwidth, height=5cm,
    ymin=0, ymax=5, xmin=0, xmax=1
]
\addplot[line width=1, color=violet]{(x^(-2)-1)^(1/2)};
\end{axis}
\end{tikzpicture}
}
\caption{The absolute value of the Leverett function for various capillary pressure models.}
\label{fig:Lev}
\end{figure}

In the pressure‐based $(\psi)$ formulation, the equations remain valid for fully saturated media ($S = 1$) but become difficult to solve as $S \to 0$. Under the Leverett relation (Eq.~\eqref{eq:leverett} and Fig.~\ref{fig:Lev}), the capillary pressure $\psi$ diverges toward $-\infty$, which makes the numerical solution unstable in very dry conditions.

The saturation‐based $(S)$ formulation avoids this difficulty at low $S$ and remains well‐defined in the dry limit ($S = 0$). However, for constitutive relationships such as van Genuchten and Haverkamp models, the capillary curve $J(S)$ is not differentiable at $S = 1$, which makes it hard to solve numerically without any regularization technique.

In practice, the choice of primary variable can strongly influence how easily the equations can be solved, especially when the domain includes both dry and saturated regions.

\subsection{A new formulation}

We propose a new formulation in which all terms of the Richards equation remain bounded. For all models presented in Table~\ref{tab:Leverett function}, the derivative of the Leverett function can be written as
\begin{equation*}
    J'(S) = C(\boldsymbol{x}) S^{-a} (1 - S^c)^{-b},
\end{equation*}
where $C(\boldsymbol{x}) > 0$, $a > 0$, $0 \le b < 1$ and $c \ge 1$. Define
\begin{equation*}
    u = \int_0^S (1 - s^c)^{-b} \, ds.
\end{equation*}
This expression can be rewritten using the incomplete beta function:
\begin{equation*}
    u = \frac{1}{c} \int_0^{S^c} s^{\frac{1}{c} - 1}(1 - s)^{(1 - b) - 1} \, ds = \frac{1}{c}\mathcal{B}\left(S^c, \frac{1}{c}, 1 - b\right),
\end{equation*}
where the incomplete beta function $\mathcal{B}$ is defined for $p, q > 0$ and $w \in [0, 1]$ as
\begin{equation*}
    \mathcal{B}(w, p, q) = \int_0^w s^{p-1} (1 - s)^{q - 1} \, ds.
\end{equation*}
The variable $u$ is well defined and monotone for $S \in [0, 1]$ and satisfies
\begin{equation*}
    0 \le u \le \frac{1}{c} \mathcal{B}\left(1, \frac{1}{c}, 1 - b\right)<\infty,
\end{equation*}
so the variable $u$ is always bounded. Defining $S = F^{-1}(u)$ with $F(S) = u(S)$, the Richards equation can be expressed in terms of $u$ as follows:
\begin{equation}
\begin{aligned}
\phi \frac{\partial S(u)}{\partial t} &- \nabla \cdot \left(K_s(\boldsymbol{x}) K_r(\boldsymbol{x}, S(u)) \left[ h_{\text{cap}}(\boldsymbol{x}) C(\boldsymbol{x}) S(u)^{-a} \left( \nabla u + \frac{\partial F}{\partial \boldsymbol{x}}(\boldsymbol{x}, S(u)) \right) \right. \right. \\
&\left. \left. + \frac{\partial \Psi}{\partial \boldsymbol{x}}(\boldsymbol{x}, S(u)) + \boldsymbol{e}_z \right] \right) = 0.
\end{aligned}
\label{eq:hetero}
\end{equation}
We have numerically solved equation \eqref{eq:hetero}; however, the terms $\frac{\partial F}{\partial \boldsymbol{x}}$ and $\frac{\partial\Psi}{\partial \boldsymbol{x}}$ introduce Dirac masses when the soil parameters are discontinuous (e.g. layered soil), complicating the application of the continuous finite element method. To address this issue, we employ a domain decomposition method instead.
For a homogeneous medium, the formulation simplifies to
\begin{equation}
\phi \frac{\partial S(u)}{\partial t} - \nabla \cdot \left(K_s K_r(S(u)) \left(h_{\text{cap}} C S(u)^{-a} \nabla u + \boldsymbol{e}_z \right)\right) = 0.
\label{eq:richarduform}
\end{equation}
To ensure that all terms in Equation~\eqref{eq:richarduform} remain bounded, it is sufficient that
\begin{equation*}
    0 \;\leq\; \lim_{S \to 0} K(S)\,S^{-a} \;<\; \infty.
\end{equation*}
This condition is satisfied by:
\begin{itemize}
    \item the Gardner model,
    \item the van Genuchten model, since \(n > 1\),
    \item the Brooks--Corey and power-law models provided \(B\geq \frac{1}{\lambda} + 1 \),
    \item the Haverkamp model provided \(\gamma \ge \beta + 1\),
    \item the van Genuchten model for capillary pressure together with the power-law model for relative permeability, provided \(B \geq \tfrac{1}{m}\).
\end{itemize}

\begin{figure}[H]
\centering
\begin{tikzpicture}
\begin{axis}[
    width=10cm,
    height=7cm, 
    grid=major,
    legend pos=south east,
    xlabel={$S$},
    ylabel={$u(S)$}
]
\addplot[blue, thick]
table [x=s, y=y, col sep=space] {integral_representation_data.txt};
\addlegendentry{$m \mathcal{B}(S^{1/m}, m, 1/n)$}
\end{axis}
\end{tikzpicture}
\caption{Graph of $u$ as a function of $S$.}
\label{fig:u for van}
\end{figure}
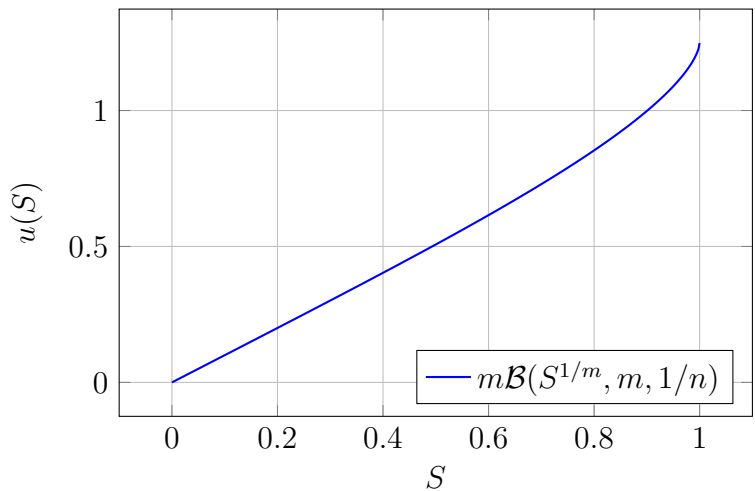

Table~\ref{tab:u variable} shows expressions of $u(S)$ for different models and Figure~\ref{fig:u for van} displays the graph of $u$ versus $S$ for the van Genuchten model.

\begin{table}[H]
\centering
\caption{$u$ as a function of $S$ for different models.}
\label{tab:u variable}
\begin{tabular}{|c|c|c|}
\hline
Model & $J'(S)$ & $u(S)$ \\
\hline
Gardner & $S^{-1}$ & $S$ \\
\hline
Brooks-Corey & $\frac{1}{\lambda} S^{-\frac{1}{\lambda} - 1}$ & $S$ \\
\hline
Haverkamp & $\frac{1}{\beta} S^{-(1 + 1/\beta)} (1 - S)^{1/\beta - 1}$ & $\beta \left[1 - (1 - S)^{1/\beta} \right]$ \\
\hline
van Genuchten & $\frac{1}{nm} S^{-1/m} (1 - S^{1/m})^{-m}$ & $m \mathcal{B}(S^{1/m}, m, 1/n)$ \\
\hline
\end{tabular}
\end{table}

Equation~\eqref{eq:richards} is supplemented by the initial condition
\begin{equation*}
    S(\boldsymbol{x}, 0) = S_0(\boldsymbol{x}), \quad \text{in } \Omega,
\end{equation*}
and boundary conditions:
\begin{equation}
\left\{
\begin{aligned}
    &S = S_D, && \text{on } \Gamma_D \times I, \\
    &\boldsymbol{q} \cdot \boldsymbol{n} = q_N, && \text{on } \Gamma_N \times I, \\
    &\nabla \Psi \cdot \boldsymbol{n} = 0, && \text{on } \Gamma_F \times I,
\end{aligned}
\right.
\label{bcrichard}
\end{equation}
where $\Gamma_D$, $\Gamma_N$ and $\Gamma_F$ are the Dirichlet, Neumann and free drainage boundaries and $S_D$, $q_N$ are prescribed data.

The Dirichlet condition on $\Gamma_D$ specifies saturation or pressure, while the Neumann condition on $\Gamma_N$ prescribes the flux. The free drainage condition on $\Gamma_F$ assumes zero normal pressure gradient.

\begin{figure}[H]
\centering
\begin{tikzpicture}[scale=0.8]
\draw[<-, dashed] (0,2) node[left] {$z$} -- (0,-9) ;
\draw[->, dashed] (-2,0) -- (10,0) node[above] {$x$}; 
\draw[color=blue] (0,0) -- (8,0);
\draw[->, color=blue, thick,dashed] (2,1) -- (2,0.5);
\draw[->, color=blue, thick, dashed] (1,1) -- (1,0.5);
\draw[->, color=blue, thick, dashed] (6,1) -- (6,0.5);
\draw[->, color=blue, thick, dashed] (7,1) -- (7,0.5);
\node[color=blue] at (4,1.5) {Fixed water flux $q_N$ or fixed saturation $S_D$};
\node[color=blue] at (4,-0.5) {$\Gamma_N$ or $\Gamma_D$};
\draw[color=red] (8,0) -- (8,-8);
\foreach \y in {-2,-1,-6,-7} {
    \draw[->, color=red, thick] (8.25,\y) -- (9,\y);
    \draw[color=red, thick] (8.5,\y-0.25) -- (8.75,\y+0.25);
}
\node[rotate=-90, color=red] at (9.5, -4) {No water flux a.e. $q_N=0$};
\node[color=red] at (7.7, -4) {$\Gamma_N$};
\draw[color=green] (8,-8) -- (0,-8);
\foreach \x in {1,2,6,7} {
    \draw[<-, color=green, thick] (\x,-9) -- (\x,-8.5);
}
\node[color=green] at (4,-9.5) {Free drainage: $\nabla \Psi \cdot \boldsymbol{n}= 0$};
\node[color=green] at (4,-7.7) {$\Gamma_F$};
\draw[color=red] (0,0) -- (0,-8);
\foreach \y in {-2,-1,-6,-7} {
    \draw[->, color=red, thick] (-0.25,\y) -- (-1,\y);
    \draw[color=red, thick] (-0.5,\y-0.25) -- (-0.75,\y+0.25);
}
\node[rotate=90, color=red] at (-1.5, -4) {No water flux a.e. $q_N=0$};
\node[color=red] at (0.5, -4) {$\Gamma_N$};
\node at (4,-4) {$\Omega$};
\node at (4,-3) {$S = S_0$ at $t = 0$};
\end{tikzpicture}
\caption{The boundary conditions.}
\label{fig:bcrich}
\end{figure}
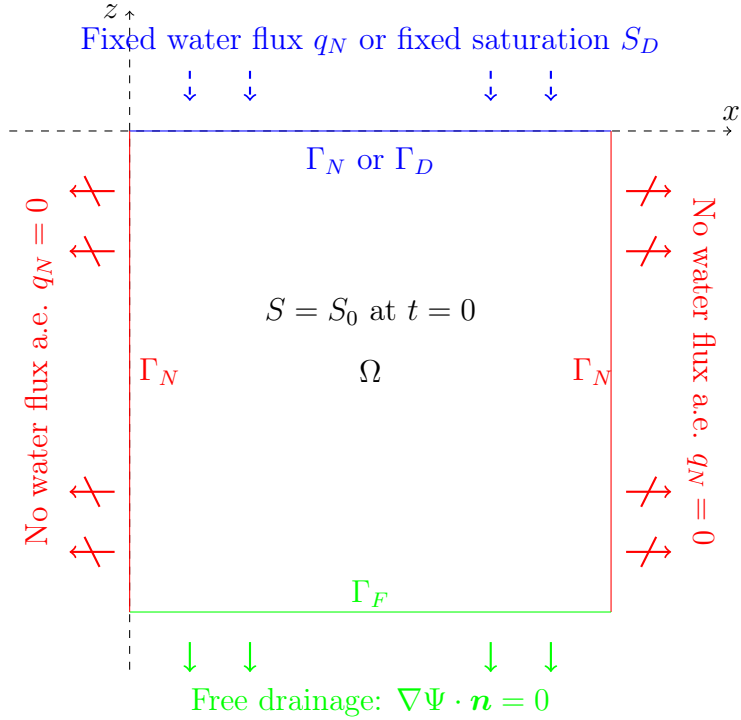

Figure~\ref{fig:bcrich} illustrates the application of these boundary conditions.

\subsection{Domain decomposition}

\begin{figure}[H]
\centering
\begin{tikzpicture}[scale=1]
        \draw[<-, dashed] (0,2) node[left] {$z$} -- (0,-9);
        \draw[->, dashed] (-2,0) -- (10,0) node[above] {$x$}; 
        \draw[color=blue] (0,0) -- (8,0);
        \draw[->, color=blue, thick, dashed] (2,1) -- (2,0.5);
        \draw[->, color=blue, thick, dashed] (1,1) -- (1,0.5);
        \draw[->, color=blue, thick, dashed] (6,1) -- (6,0.5);
        \draw[->, color=blue, thick, dashed] (7,1) -- (7,0.5);
        \node[color=blue] at (4,1.5) {Fixed water flux $q_N$ or a fixed water saturation $S_D$};
        \node[color=blue] at (4,-0.5) {$\Gamma_N$ or $\Gamma_D$};
        \draw[color=red]  (8,0) -- (8,-4);
        \draw[->, color=red, thick] (8.25,-2) -- (9,-2);
        \draw[color=red, thick] (8.5,-2.25) -- (8.75,-1.75);
        \node[rotate=-90, color=red] at (9.5,-4) {No water flux a.e.\ $q_N=0$};
        \node[color=red] at (7.7,-2) {$\Gamma_N$};

        \draw[color=red] (0,0) -- (0,-4);
        \draw[->, color=red, thick] (-0.25,-2) -- (-1,-2);
        \draw[color=red, thick] (-0.5,-2.25) -- (-0.75,-1.75);
        \node[color=red] at (0.5,-2) {$\Gamma_N$};
        \fill[blue, opacity=0.3, domain=0:8] (8,0) -- (0,0) -- plot ({\x},{0.5*sin(\x r)-4});
        \node at (4,-1.5) {$u_1(\boldsymbol{x}, t)$};
        \node at (4,-2) {$\Omega_1$};
        \node at (4,-3) {${\theta_s}_1,\; {\theta_r}_1,\; {K_s}_1,\; {h_{cap}}_1, \dots$};

        \draw[domain=0:8, smooth, variable=\x, purple] plot ({\x},{0.5*sin(\x r)-4});
        \draw[color=red] (8,-4) -- (8,-8);
        \draw[->, color=red, thick] (8.25,-6) -- (9,-6);
        \draw[color=red, thick] (8.5,-6.25) -- (8.75,-5.75);
        \node[color=red] at (7.7,-6) {$\Gamma_N$};
        \draw[color=green] (8,-8) -- (0,-8);
        \draw[<-, color=green, thick] (2,-9) -- (2,-8.5);
        \draw[<-, color=green, thick] (1,-9) -- (1,-8.5);
        \draw[<-, color=green, thick] (6,-9) -- (6,-8.5);
        \draw[<-, color=green, thick] (7,-9) -- (7,-8.5);
        \node[color=green] at (4,-9.5) {Free drainage: \ $\nabla \Psi \cdot \boldsymbol{n_2}= 0$};
        \node[color=green] at (4,-7.7) {$\Gamma_F$ or $\Gamma_D$};
        \draw[color=red] (0,-4) -- (0,-8);
        \draw[->, color=red, thick] (-0.25,-6) -- (-1,-6);
        \draw[color=red, thick] (-0.5,-6.25) -- (-0.75,-5.75);
        \node[rotate=90, color=red] at (-1.5,-4) {No water flux a.e.\ $q_N=0$};
        \node[color=red] at (0.5,-6) {$\Gamma_N$};
        \fill[yellow, opacity=0.3, domain=0:8] (8,-8) -- (0,-8) -- plot ({\x},{0.5*sin(\x r)-4});
        \node at (4,-5.5) {$u_2(\boldsymbol{x}, t)$};
        \node at (4,-6) {$\Omega_2$};
        \node at (4,-7) {${\theta_s}_2,\; {\theta_r}_2,\; {K_s}_2,\; {h_{cap}}_2,\dots$};

        \node at (4,-4.75) {$\Gamma$};
    \end{tikzpicture}
\caption{Layered soil problem.}
\label{fig:layeredsoil}
\end{figure}
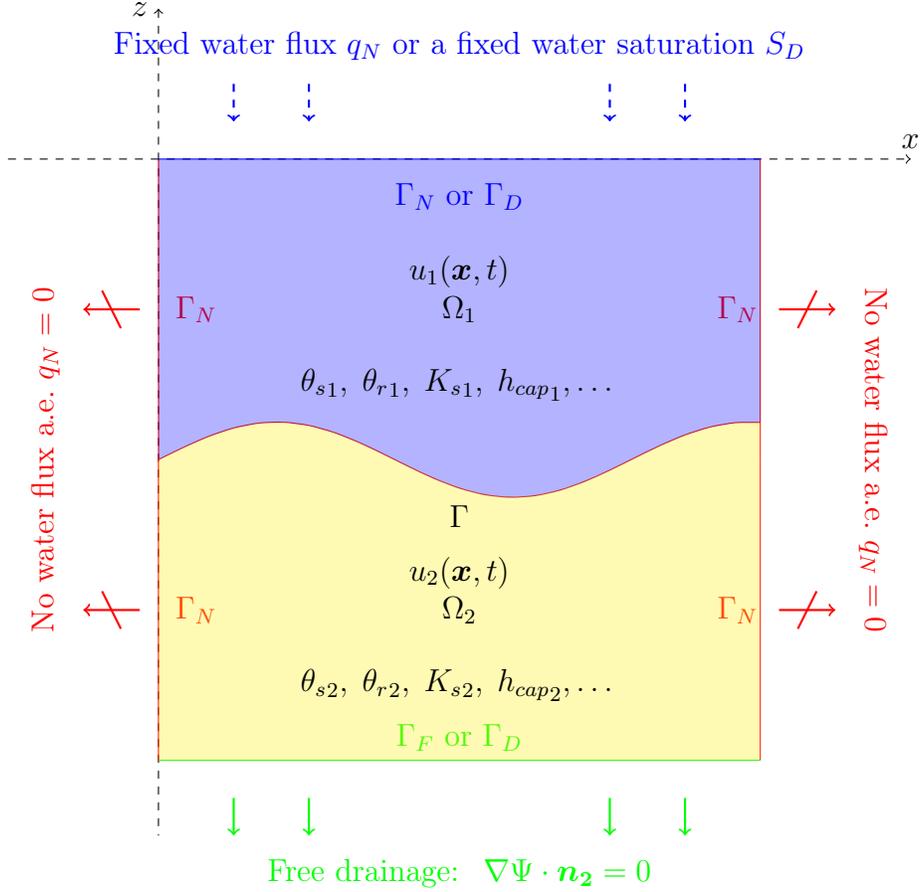

Layered soil problems can be efficiently solved using domain decomposition methods. Let the domain be $\overline{\Omega} = \overline{\Omega_1 \cup \Omega_2}$, where $\Omega_1 \cap \Omega_2 = \varnothing$ and $\partial \Omega_1 \cap \partial \Omega_2 = \Gamma$ (see Figure~\ref{fig:layeredsoil}). The governing equations can be solved independently in subdomains $\Omega_1$ and $\Omega_2$ and coupled across the interface $\Gamma$.

Figure~\ref{fig:layeredsoil} shows a vertical cross-section of a layered soil profile with two distinct layers, $\Omega_1$ (top) and $\Omega_2$ (bottom), separated by the interface $\Gamma$. The top boundary has a specified water flux or saturation, the lateral boundaries are impermeable and the bottom boundary allows for free drainage. Each layer has distinct state variable $u_i$, $i=1,2$, and soil-specific hydraulic properties. The arrows indicate the direction of flow and the nature of the boundary conditions.

We consider a domain decomposition for the Richards equation using the van Genuchten model. 
Assuming continuity of $\Psi$, the coupled problem reads:
\begin{equation}
\left\{
\begin{aligned}
    &\phi_i \frac{\partial S_i(u_i)}{\partial t} - \nabla \cdot \left(K_{s,i} K_{r,i}(S_i(u_i)) \left(h_{\text{cap},i} C_i S_i(u_i)^{-a} \nabla u_i + \boldsymbol{e}_z \right)\right) = 0, && \text{in } \Omega_i \times I, \; i=1,2, \\
    &u_i(\boldsymbol{x}, 0) = u_{i,0}(\boldsymbol{x}), && \text{in } \Omega_i, \\
    &u_1 = u_D, && \text{on } \Gamma_D \times I, \\
    &\boldsymbol{q}_1 \cdot \boldsymbol{n}_1 = q_N, && \text{on } \Gamma_N \cap \partial \Omega_1 \times I, \\
    &\boldsymbol{q}_2 \cdot \boldsymbol{n}_2 = q_N, && \text{on } \Gamma_N \cap \partial \Omega_2 \times I, \\
    &\nabla \Psi_2 \cdot \boldsymbol{n}_2 = 0, && \text{on } \Gamma_F \times I, \\
    &\Psi_1(S_1(u_1)) = \Psi_2(S_2(u_2)), \qquad 
     \boldsymbol{q}_1 \cdot \boldsymbol{n}_1 = \boldsymbol{q}_2 \cdot \boldsymbol{n}_2, 
     && \text{on } \Gamma \times I.
\end{aligned}
\right.
\label{eq:domaindecomp}
\end{equation}

On each subdomain, we thus solve the same Richards equation 
with the boundary and initial conditions above, 
but the interface conditions on $\Gamma \times I$ can be equivalently written as
\begin{equation}
\left\{
\begin{aligned}
    &u_1 = S_1^{-1}\!\left(\Psi_1^{-1}\!\big(\Psi_2(S_2(u_2))\big)\right), \\
    &u_2 = S_2^{-1}\!\left(\Psi_2^{-1}\!\big(\Psi_1(S_1(u_1))\big)\right), \\
    &\boldsymbol{q}_1 \cdot \boldsymbol{n}_1 = \boldsymbol{q}_2 \cdot \boldsymbol{n}_2.
\end{aligned}
\right.
\label{eq:interface}
\end{equation}

\subsection{Variational formulations}
In this section, we present the variational formulations for the \(\Psi\), \(S\) and \(u\) formulations.

\subsubsection{\(\Psi\)-formulation}
The \(\Psi\)-formulation is given by: find \(\Psi \in L^2(I; H^1(\Omega))\) such that
\begin{equation*}
\left\{
    \begin{aligned}
        & \Psi - \Psi_D \in L^2(I; H^1_{\Gamma_D}(\Omega)), \\
        & \int_\Omega \phi \frac{\partial S(\Psi)}{\partial t} v\,dx + \int_\Omega K_s K_r(\Psi)\nabla (\Psi + z)\cdot\nabla v\,dx \\
        & \quad + \int_{\Gamma_N} q_N v\,ds + \int_{\Gamma_F} K_s K_r(\Psi)\boldsymbol{e}_z\cdot\boldsymbol{n} v\,ds = 0,\quad\forall v\in H^1_{\Gamma_D}(\Omega) \text{ a.e. in } I.
    \end{aligned}
\right.
\end{equation*}

\subsubsection*{\(S\)-formulation}
The \(S\)-formulation, in the case of a homogeneous medium, is given by: find \(S \in L^2(I; H^1(\Omega))\) such that
\begin{equation}
\left\{
    \begin{aligned}
        & S - S_D \in L^2(I; H^1_{\Gamma_D}(\Omega)), \\
        & \int_\Omega \phi \frac{\partial S}{\partial t} v\,dx + \int_\Omega K_s K_r(S)\left(h_{cap}J'(S)\nabla S + \boldsymbol{e}_z\right)\cdot\nabla v\,dx \\
        & \quad + \int_{\Gamma_N} q_N v\,ds + \int_{\Gamma_F} K_s K_r(S)\boldsymbol{e}_z\cdot\boldsymbol{n} v\,ds = 0,\quad\forall v\in H^1_{\Gamma_D}(\Omega)  \text{ a.e. in } I.
    \end{aligned}
\right.
\label{var:Sform}
\end{equation}

\subsubsection*{\(u\)-formulation}
The variational formulation in terms of \(u\), in the case of a homogeneous medium, is given by: find \(u\in L^2(I; H^1(\Omega))\), such that
\begin{equation}
\left\{
    \begin{aligned}
        & u - u_D \in L^2(I; H^1_{\Gamma_D}(\Omega)), \\
        & \int_\Omega \phi \frac{\partial S(u)}{\partial t} v\,dx + \int_\Omega K_s K_r(S(u))\left(h_{cap}S(u)^{-a}\nabla u + \boldsymbol{e}_z\right)\cdot\nabla v\,dx \\
        & \quad + \int_{\Gamma_N} q_N v\,ds + \int_{\Gamma_F} K_s K_r(S(u))\boldsymbol{e}_z\cdot\boldsymbol{n} v\,ds = 0,\quad\forall v\in H^1_{\Gamma_D}(\Omega)  \text{ a.e. in } I.
    \end{aligned}
\right.
\label{eq:uhomoformulation}
\end{equation}
In all formulations, we require that \(\phi \, \frac{\partial S}{\partial t} \in L^2(I; L^2(\Omega))\).

\section{Numerical methods}
\label{sec:NMETHODS}
\subsubsection*{Temporal discretization}
Let $N > 0$. The interval $I = (0,T)$ is subdivided into $N$ sub-intervals,
\[
t_n = n\Delta t, \quad \Delta t = \frac{T}{N},
\]
and we define $S^n \approx S(\cdot, t_n)$, $\Psi^n \approx \Psi(\cdot, t_n)$ and $u^n=u(\cdot,t_n)$.

\subsubsection*{Spatial discretization}
It is assumed that $\overline{\Omega}$ is a polygonal domain and $\Gamma \subset \Omega$ a curve linear by parts. We discretize the domain $\overline{\Omega}$ by partitioning it into elements $\Omega^e$. We denote by $\tau_h$ the set of all elements $\Omega^e$. We suppose that $\tau_h$ is a conformal mesh \cite{ErnGuermond2021}. The choice of the elements $\Omega^e$ depends on the space dimension $d$. If $d=1$, the $\Omega^e$ are segments, if $d=2$, then they are triangles and for $d>2$, they are closed $d$-simplices. We define the mesh size $h$ by
\[
h = \underset{\tau_h}{\max} \, \mathrm{diam}(\Omega^e),
\]
where $\mathrm{diam}$ denotes the diameter of the element $\Omega^e$ given by
\[
\mathrm{diam}(\Omega^e) = \underset{x,y \in \Omega^e}{\sup} \lVert x-y \rVert _2.
\]
The Galerkin finite element space $V_h(\Omega) $ is defined as follows:
\[
V_h(\Omega) = \left\{ v_h \in H^1(\Omega) \mid v_{h|\Omega^e} \in \mathcal{P}_k(\Omega^e), \forall \Omega^e \in \mathcal{T}_h \right\},
\]
and if $\Gamma \subset \partial\Omega$ we define 
\[
V_h(\Gamma) = \{ {v_h}_{|\Gamma}|\; v_h \in V_h(\Omega) \},
\]
where $\mathcal{P}_k(\Omega^e)$ represents the space of polynomials of degree $k$ on the element $\Omega^e$. Here, we use linear elements $k=1$. We approximate $\Psi$, $S$ and  $u$ by $\Psi_h^n$, $S_h^n$ $u_h^n \in V_h$, respectively. If there is only one domain, we denote $V_h = V_h(\Omega)$.
\subsubsection*{Time stepping schemes for the Richards equation}
In the following section, we give time and space discretizations for the $S$ based formulation and the $u$ based formulation.

\subsubsection*{$S$-formulation}
Using the finite element method, we discretize the variational formulation \eqref{var:Sform}. For the temporal discretization, we use a semi-implicit Euler method. Given $S^n_h\in V_h$, find $S^{n+1}_h\in V_h$ such that,
\begin{equation*}
\left\{
    \begin{aligned}
        &S_h^0 = \Pi_h S_0,\\
        &\int_\Omega \phi \frac{S^{n+1}_h-S^{n}_h}{\Delta t} v_hdx+ \int_\Omega K_sK_r(S^{n}_h)\left(h_{cap}J'(S^n_h)\nabla S^{n+1}_h +\boldsymbol{e}_z\right) \cdot \nabla v_h dx \\
        &+\int_{\Gamma_N} q_N v_h ds + \int_{\Gamma_F} K_s K_r(S^{n}_h) \boldsymbol{e}_z \cdot \boldsymbol{n} v_hds  =0,\; \forall v_h \in V_h.
    \end{aligned}
    \right.
\end{equation*}
where $\Pi_h: H^1(\Omega) \to V_h$ is a projection operator onto $V_h$. 
It is also possible to use an implicit Euler method, where all expressions are evaluated at $t_{n+1}$.

\subsubsection*{$u$-formulation}
We discretize the variational formulation \eqref{eq:uhomoformulation}, using the finite element method. For the temporal discretization, we employ a semi-implicit Euler scheme. Newton's method to solve the resulting nonlinear system. Let $u^n_h\in V_h$ be given for $n=0, \dots, N-1$, find $u^{n+1}_h\in V_h$ such that
\begin{equation}
\left\{
    \begin{aligned}
        & u_0^h = \Pi_h u_0,\\
        &\int_\Omega \phi\frac{S(u^{n+1}_h)-S(u^n_h)}{\Delta t} v_hdx \\
        &+ \int_\Omega K_sK_r(S(u^{n}_h))\left(h_{cap}S(u^{n}_h)^{-a}\nabla u^{n+1}_h +\boldsymbol{e}_z\right) \cdot \nabla v_h dx \\
        &+\int_{\Gamma_N} q_N v_h ds + \int_{\Gamma_F} K_s K_r(S(u^{n}_h)) \boldsymbol{e}_z \cdot \boldsymbol{n} v_hds  =0,\; \forall v_h \in V_h,
    \end{aligned}
    \right.
    \label{eq: disuhomoformulation}
\end{equation}
To solve the nonlinear system \eqref{eq: disuhomoformulation}, we will use Newton's method. Given the initial guess $u^{n+1,0}_h = u^n_h$, assume $u^{n+1,m}_h \in V_h$ is known for $m\ge 0$, find $u^{n+1,m+1}_h \in V_h$ such that
\begin{equation}
\left\{
    \begin{aligned}
        &\int_\Omega \phi\frac{S(u^{n+1, m+1}_h)-S(u^n_h)}{\Delta t} v_hdx\\
        &+ \int_\Omega K_sK_r(S(u^{n}_h))\left(h_{cap}S(u^{n}_h)^{-a}\nabla u^{n+1,m+1}_h +\boldsymbol{e}_z\right) \cdot \nabla v_h dx \\
        &+\int_{\Gamma_N} q_N v_h ds + \int_{\Gamma_F} K_s K_r(S(u^{n}_h) )\boldsymbol{e}_z \cdot \boldsymbol{n} v_hds  =0,\; \forall v_h \in V_h.
    \end{aligned}
    \right.
    \label{eq: picarduhomoformulation}
\end{equation}
We approximate $S(u^{n+1,m+1}_h)$ by using a first order Taylor expansion in $S$,
\begin{equation*}
    S(u^{n+1, m+1}_h) = S(u^{n+1, m}_h)+\frac{\partial S}{\partial u}(u^{n+1, m}_h)(u^{n+1, m+1}_h-u^{n+1, m}_h) + O((u^{n+1, m+1}_h-u^{n+1, m}_h)^2).
\end{equation*}
Thus, the nonlinear semi-implicit scheme is given by:
find $u^{n+1,m+1}_h \in V_h$ such that
\begin{equation}
\left\{
    \begin{aligned}
        &\int_\Omega \phi \frac{\partial S}{\partial u}(u^{n+1, m}_h) \frac{u^{n+1, m+1}_h-u^{n+1, m}_h}{\Delta t} v_hdx +\int_\Omega \phi \frac{S(u^{n+1, m}_h)-S(u^{n}_h)}{\Delta t} v_hdx\\ &+ \int_\Omega K_sK_r(S(u^{n}_h))\left(h_{cap}S(u^{n}_h)^{-a}\nabla u^{n+1,m+1}_h +\boldsymbol{e}_z\right) \cdot \nabla v_h dx \\
        &+\int_{\Gamma_N} q_N v_h ds + \int_{\Gamma_F} K_s K_r(S(u^{n}_h)) \boldsymbol{e}_z \cdot \boldsymbol{n} v_hds  =0,\; \forall v_h \in V_h.
    \end{aligned}
    \right.
    \label{eq:modpicarduhomoformulation}
\end{equation}
By only doing one iteration, we get a linear semi-implicit method, which is equivalent to finding $u^{n+1}_h\in V_h$ such that
\begin{equation}
\left\{
    \begin{aligned}
        &\int_\Omega \phi \frac{\partial S}{\partial u}(u^{n}_h) \frac{u^{n+1}_h-u^{n}_h}{\Delta t} v_hdx+ \int_\Omega K_sK_r(S(u^{n}_h))\left(h_{cap}S(u^{n}_h)^{-a}\nabla u^{n+1}_h +\boldsymbol{e}_z\right) \cdot \nabla v_h dx \\
        &+\int_{\Gamma_N} q_N v_h ds + \int_{\Gamma_F} K_s K_r(S(u^{n}_h)) \boldsymbol{e}_z \cdot \boldsymbol{n} v_hds  =0,\; \forall v_h \in V_h.
    \end{aligned}
    \right.
    \label{eq:semiimplicitscheme}
\end{equation}

\subsubsection*{Domain decomposition method Heterogeneous problems}
We discretize \eqref{eq:domaindecomp} in space using the finite element method on a conforming mesh, assuming that $\Gamma$ is piecewise linear. 
Time discretization is performed with a semi-implicit Euler scheme. 
The fully discrete system is nonlinear, and Newton’s method is employed to solve it.

\medskip\noindent
Given $(u_{1h}^n,u_{2h}^n)\in V_h(\Omega_1)\times V_h(\Omega_2)$, 
we compute $(u_{1h}^{n+1},u_{2h}^{n+1})\in V_h(\Omega_1)\times V_h(\Omega_2)$ such that, 
for each $i=1,2$,
\begin{equation*}
\left\{
\begin{aligned}
  &u_{ih}^0 = \Pi_h u_{i,0},\\
  &\int_{\Omega_i}\phi_i\frac{S_i(u_{ih}^{n+1})-S_i(u_{ih}^n)}{\Delta t}v_h\,dx
   +\int_{\Omega_i} K_{si}K_{ri}(S_i(u_{ih}^n))\left(h_{\text{cap},i}S_i(u_{ih}^n)^{-a}\nabla u_{ih}^{n+1}+\boldsymbol{e}_z\right)\cdot\nabla v_h\,dx \\
  &\quad +\int_{\Gamma_N}q_N v_h\,ds
   +\int_{\Gamma_F}K_{si}K_{ri}(S_i(u_{ih}^n))\boldsymbol{e}_z\cdot\boldsymbol{n}_i v_h\,ds
   \\
   &+\int_{\Gamma}\boldsymbol{q}_i^{n+1}\cdot\boldsymbol{n}_i v_h\,ds=0, \forall v_h\in V_h(\Omega_i).
\end{aligned}
\right.
\end{equation*}

\medskip\noindent
Newton’s method is applied as follows: given $u_{ih}^{n+1,m}$, 
find $u_{ih}^{n+1,m+1}\in V_h(\Omega_i)$ such that
\begin{equation*}
\left\{
\begin{aligned}
 &\int_{\Omega_i}\phi_i \frac{\partial S_i}{\partial u}(u_{ih}^{n+1,m})
     \frac{u_{ih}^{n+1,m+1}-u_{ih}^{n+1,m}}{\Delta t} v_h\,dx
 +\int_{\Omega_i}\phi_i\frac{S_i(u_{ih}^{n+1,m})-S_i(u_{ih}^n)}{\Delta t} v_h\,dx \\
 &\quad +\int_{\Omega_i} K_{si}K_{ri}(S_i(u_{ih}^n))
     \left(h_{\text{cap},i}S_i(u_{ih}^n)^{-a}\nabla u_{ih}^{n+1,m+1}+\boldsymbol{e}_z\right)\cdot\nabla v_h\,dx \\
 &\quad +\int_{\Gamma_N} q_N v_h\,ds
 +\int_{\Gamma_F} K_{si}K_{ri}(S_i(u_{ih}^n))\boldsymbol{e}_z\cdot\boldsymbol{n}_i v_h\,ds
 \\
 &+\int_{\Gamma}\boldsymbol{q}_i^{n+1,m+1}\cdot\boldsymbol{n}_i v_h\,ds = 0,  \forall v_h\in V_h(\Omega_i).
\end{aligned}
\right.
\end{equation*}

\medskip\noindent
To enforce the interface coupling, we introduce Robin-type transmission conditions. 
At iteration $m+1$, the problem on $\Omega_1$ is solved using
\begin{equation*}
\boldsymbol{q}_{1}^{n+1,m+1}\cdot\boldsymbol{n}_{1}+\lambda u^{n+1,m+1}_{1h} 
= -\boldsymbol{q}_{2}^{n+1,m}\cdot\boldsymbol{n}_{2}
+ \lambda S_1^{-1}\!\left(\Psi_1^{-1}\!\big(\Psi_2(S_2(u^{n+1,m}_{2h}))\big)\right),
\qquad \text{on }\Gamma,
\end{equation*}
which yields $u^{n+1,m+1}_{1h}$.  
Then, the problem on $\Omega_2$ is solved using
\begin{equation*}
\boldsymbol{q}_{2}^{n+1,m+1}\cdot\boldsymbol{n}_{2}+\lambda u^{n+1,m+1}_{2h} 
= -\boldsymbol{q}_{1}^{n+1,m+1}\cdot\boldsymbol{n}_{1}
+ \lambda S_2^{-1}\!\left(\Psi_2^{-1}\!\big(\Psi_1(S_1(u^{n+1,m+1}_{1h}))\big)\right),
\qquad \text{on }\Gamma,
\end{equation*}

\section{Numerical tests}
\label{sec:Numres}
In this section, we present a series of numerical test cases to evaluate the performance of the proposed finite element method. The objective is to demonstrate the capability of the $u$-formulation to solve problems under fully saturated and dry conditions. Each test case focuses on a specific scenario within this framework.

The implementation was performed in FreeFEM++~\cite{Hecht2012},  or FENICSx~\cite{Dolfinx23} with Paraview and Python (Matplotlib) used for visualization. The resulting linear systems were solved using the sparse direct solvers UMFPACK~\cite{Davis2004}, PETSc and MUMPS~\cite{Amestoy2001, Amestoy2019, Balay2023}.
\subsection*{Order of convergence analysis}
To compute the error between the numerical and reference solutions, we used the \(L^2\) errors at the final time $T$:
\[
E_{h,\Delta t}(T)^2=\|u(\cdot,T) - u_{h,\Delta t}(\cdot,T)\|_{L^2(\Omega)}^2 = \int_\Omega |u(\boldsymbol{x},T) - u_{h,\Delta t}(\boldsymbol{x},T)|^2d\boldsymbol{x},
\]
where \(u\) represents the reference solution, and \(u_{h,\Delta t}\) denotes the numerical approximation computed with mesh size $h$ and time step $\Delta t$. 
To compute the order of the method in space, $p$, we fix a small time step, $\Delta t$, and calculate $p$ as follows:
$$
p=\log_2\left(E_{h,\Delta t}(T)/E_{h/2,\Delta t}(T)\right),
$$
and the order in time, after fixing a small space step $h$, is computed by:
$$
q=\log_2\left(E_{h,\Delta t}(T)/E_{h,\Delta t/2}(T)\right).
$$

\subsection{1D fully unsaturated fibrous sheets}
\label{subsec:1Dfullyunsaturated}
This test case evaluates the performance of the $u$-formulation in a dry setting. The goal is to demonstrate that the method is stable and accurate. Results obtained with the $u$-formulation are compared against those from the classical $S$-formulation to assess consistency and reliability.

We follow the setup in \cite{Jaganathan2009} for modeling mineral oil flow in fibrous sheets. We employ the \(S\)-formulation in one dimension:
\[
\phi\frac{\partial S}{\partial t} - \frac{\partial}{\partial z}\left(K_sK_r(S)h_{cap}J'(S)\frac{\partial S}{\partial z} \right) - \rho g K_s \frac{\partial K_r(S)}{\partial z} = 0.
\]
We then apply our \(u\)-formulation to validate its effectiveness in fully dry domains:
\[
\phi\frac{\partial S(u)}{\partial t} - \frac{\partial}{\partial z}\left(K_s h_{cap}\frac{1}{nm}K_r(S(u))S(u)^{-1/m}\frac{\partial u}{\partial z} \right) - \rho g K_s\frac{\partial K_r(S(u))}{\partial z} = 0.
\]
We use the van Genuchten model for the Leverett function and a power law for the relative permeability (see Tables~\ref{tab:Leverett function} and~\ref{tab:u variable}). The initial and boundary conditions are:
\begin{align*}
    &S(z,0) = 0, &&\forall z \in [0,100],\\
    &S(0,t) = 1, \quad S(100,t) = 0, &&\forall t \in (0, 86400].
\end{align*}
Figure~\ref{fig:1D_BC_IC} illustrates the domain and the boundary conditions:
\begin{figure}[H]
    \centering
    \begin{tikzpicture}[scale=0.67]
        \node[color = blue] at (1, 0) {$S=0$} ;
        \node[color = blue] at (-1.25, 0) {$z=100$} ;
        \draw[color=blue, fill=blue] (0,0) circle (2pt);
        \draw[color=red]  (0,0) -- (0,-8);
        \draw[->]  (4,-2) -- (4,-6);
        \draw[->, color=red, dashed]  (0,0) -- (0,1);
        \node at (4.5, -4) {$\mathrm{g}$} ;
        \draw[color=green, fill=green] (0,-8) circle (2pt);
        \node[color = green] at (1, -8) {$S=1$} ;
        \node[color = green] at (-1, -8) {$z=0$} ;
        \node[color = red] at (1, -4) {$S=S_0$} ;
    \end{tikzpicture}
    \caption{1D boundary conditions.}
    \label{fig:1D_BC_IC}
\end{figure}
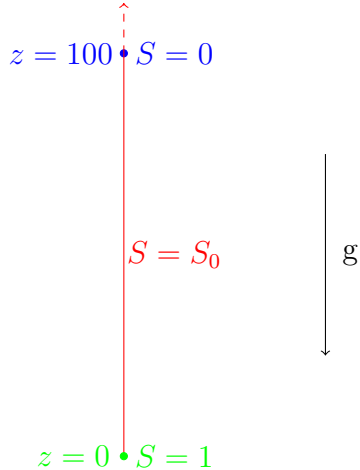

The parameters of the porous fibrous sheet used in this study are given in table~\ref{tab:solparam}, as reported in \cite{Jaganathan2009}.

\begin{table}[H]
    \centering
    \caption{Values of soil hydraulic parameters}
    \begin{tabular}{|c|c|c|c|c|c|c|c|}
         \hline
         $\theta_s$ & $\theta_r$ & $n$ & $K_s$ [mm$^3$/g/s] & $h_{cap}$ [g/mm/s$^2$] & $B$ & $\rho$ [g/mm$^2$] & $g$ [mm/s$^2$] \\
         \hline
         0.95 & 0 & 4 & 0.01 & 342 & 4.6 & $8.41\times 10^{-4}$ & $9.8\times 10^{3}$ \\
         \hline
    \end{tabular}
    \label{tab:solparam}
\end{table}

\begin{figure}[H]
    \centering
    \includegraphics[width=1\linewidth]{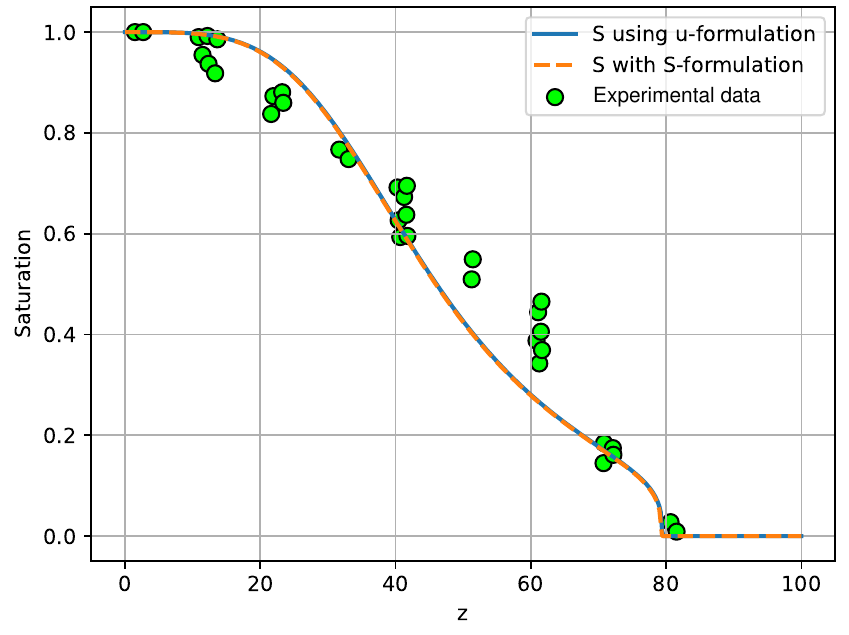}
    \caption{Comparison between the results obtained using linear semi-implicit \(u\)- and \(S\)-formulations for \(h = 1/5\)~mm and \(\Delta t = 1\) s, and with experimental data.}
    \label{fig:fibrous sheets}
\end{figure}

The reference solution based on the \(S\)-formulation in \cite{Jaganathan2009} is compared with the numerical solution of the linear semi-implicit \eqref{eq:semiimplicitscheme}. Computations are carried out using a uniform 1D mesh with 500 points and a time step \(\Delta t = 1\)s. The results presented in Figure~\ref{fig:fibrous sheets} demonstrate excellent agreement between the $u$- and $S$-formulations, and with experimental data taken from \cite{Jaganathan2009}. The linear semi-implicit scheme performs well when \(S < 1\), which holds at all internal grid points for this test case. Consequently, the S-formulation works well for this test case.

\subsection{Test case with 1D fully saturated and fully unsatured regions.}

We conduct a 1D numerical test case involving both fully saturated and fully unsaturated regions to illustrate that the $u$-formulation handles such configurations properly. None of the $S$-formulation or $\psi$-formulation can be applied to this test case unless special switching or regularization strategies are introduced.

The van Genuchten–Mualem model is used for the constitutive relationships (see Tables~\ref{tab:Leverett function} and~\ref{tab:u variable}). Initial values and boundary conditions are given as follows

\begin{align*}
    &S(z,0)    = \left\{\begin{aligned}
        &0, \quad \text{if} \; z \in [0, 0.5],\\
        &1, \quad \text{if} \; z \in (0.5, 1],
    \end{aligned} \right.
    \\&
    \\&\;S(1,t)    = 1,\;S(0, t) = 0,\; \forall t \in (0, 1].
\end{align*}
A geometry similar to that presented in Figure~\ref{fig:1D_BC_IC} was considered, and the clay loam soil parameters listed in Table~\ref{tab:solparam1} were taken for use in the Hydrus-1D software.

\begin{table}[H]
    \centering
    \caption{Values of soil hydraulic parametres }
    \begin{tabular}{|c|c|c|c|c|c|c|c|c|c|}
         \hline$\theta_s$ & $\theta_r$ & $n$ & $K_s[\mathrm{m/days}]$ & $\alpha [\mathrm{m}^{-1}]$ &$h_{cap}[\mathrm{m}]$\\
\hline 0.41 & 0.095  & 1.31 & 0.0624 & 1.9 & $\alpha^{-1}$\\
\hline
    \end{tabular}
    
    \label{tab:solparam1}
\end{table}
For the present configuration, we employ the nonlinear semi-implicit scheme \eqref{eq:modpicarduhomoformulation}. 
The initial condition is non-smooth and reaches the value $1$ on parts of the domain. 
This creates stiffness and degeneracy. 
The linear semi-implicit \eqref{eq:semiimplicitscheme} treatment would not handle this situation robustly. In Section~\ref{subsec:1Dfullyunsaturated}, the situation is different. 
There, the initial solution is smooth and does not reach the value $1$ on any subregion. 

The test case presented here uses the $u$-formulation with a uniform 1D mesh of 5000 points. 
We adopt a time step of $10^{-5}$ hours. 
The solver is set with an absolute tolerance of $10^{-8}$, a relative tolerance of $10^{-8}$, and a maximum of 100 Newton iterations. 
On average, convergence is achieved within 1--5 iterations. 

\begin{figure}[H]
    \centering
    \includegraphics[width=1\linewidth]{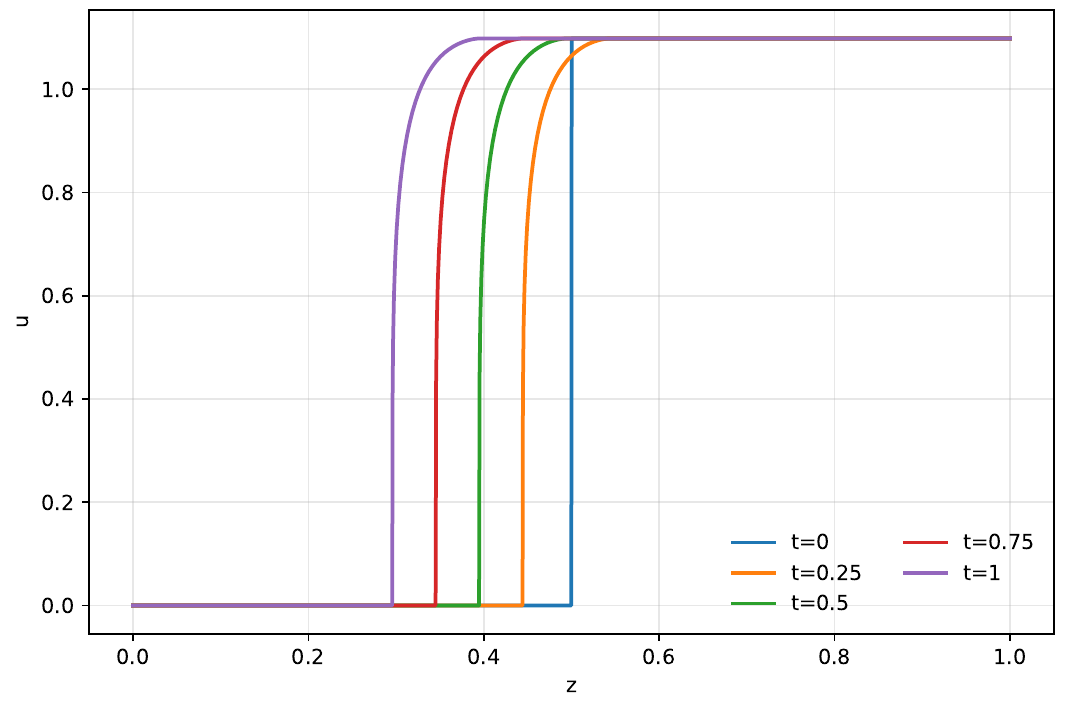}
    \caption{the solution $u$ obtained by the nonlinear semi-implicit $u$-formulation, using $h=2\times10^{-4} \text{ cm}$ and $\Delta t=10^{-5} \text{ hours}$.}
     \label{fig:results}
\end{figure}

Figure~\ref{fig:results} shows the results obtained with the implicit method. 
The solution produced by our approach is monotone. 

Computing the solution on a mesh of size $h=4\times10^{-4}$ and $dt=2\times10^{-5}$, we get an error of $3.340737\times 10^{-3}$ at the final time, when comparing with the solution from figure ~\ref{fig:results}, which shows independence of the solution from mesh and time step.
To further validate this solution, we use the semi-implicit $S$-formulation method, where we regularize the derivative of $J$ like the following, given $\delta>0$
\begin{equation*}
    J'(S) = \left\{\begin{aligned}
        &J'(1-\delta), \;S>1-\delta,\\
        &J'(S), \;S<1-\delta.
    \end{aligned} \right.
\end{equation*}
Using a very small regularization parameter $\delta=10^{-16}$, we take the saturation $S$ computed with the nonlinear semi-implicit $u$-formulation as a reference. We then verify that the solutions obtained with the regularized semi-implicit $S$-formulation converge towards this reference as the time step $\Delta t$ is refined with a fixed $h=2.0\times 10^{-4}$. 

\begin{figure}[H]
    \centering
    \includegraphics[width=1\linewidth]{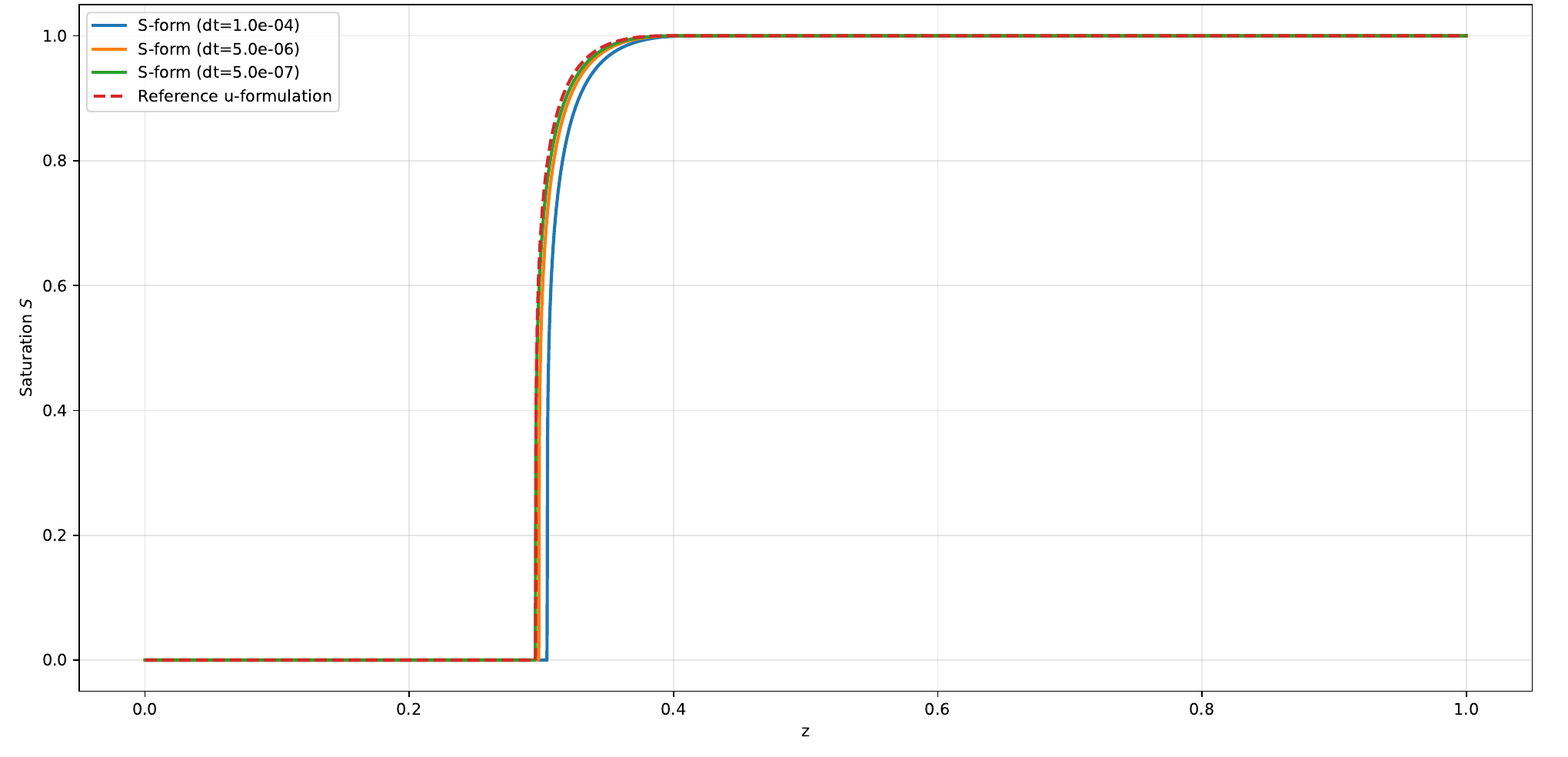}
    \caption{Comparison of $S$ computed using the $u$-formulation and the regularized $S$-formulation for different time steps at the final time.}
    \label{fig:svsu}
\end{figure}

\begin{table}[H]
\centering
\begin{tabular}{cc}
\hline
\textbf{$\Delta t$} & $\| \text{error} \|_{L^2}$ \\
\hline
$1.0 \times 10^{-4}$ & $6.808641 \times 10^{-2}$ \\
$5.0 \times 10^{-6}$ & $2.972399 \times 10^{-2}$ \\
$5.0 \times 10^{-7}$ & $5.447313 \times 10^{-3}$ \\
\hline
\end{tabular}
\caption{$L^2$ errors of the $S$-formulation for different time steps.}
\label{tab:sform_errors}
\end{table}

Figures~\ref{fig:svsu} and Table~\ref{tab:sform_errors} demonstrate that the regularized semi-implicit $S$-formulation converges towards the reference $u$-formulation, albeit at a relatively slow rate. Achieving higher accuracy requires very small time steps: for instance, computing with $\Delta t = 5.0 \times 10^{-7}$ takes approximately two hours on a single core. These results show that the nonlinear semi-implicit $u$-formulation method is accurate, robust and stable across both dry and fully saturated regimes.

\subsection{Manufactured solution}

We now consider a smoother configuration. This case complements the discontinuous initial conditions discussed in the previous subsection. The problem is defined on a one-dimensional spatial domain, $x \in [0,1]$ m. The time interval is $t \in [0,1]$ days. 

The soil hydraulic properties are described by the van Genuchten--Mualem model. The corresponding parameters are listed in Table~\ref{tab:soilparams_mms}. They are the same as those of clay loam, except that the value of $n$ is increased. This adjustment is made to simplify the computation of the source term. For small $n$, the relative permeability shows a strong singularity at $S=1$. In that case, many quadrature points would be needed. For larger $n$, the singularity is weaker, and the computation becomes easier.

\begin{table}[h!]
\centering
\caption{Values of soil hydraulic parameters used in the manufactured solution test case.}
\label{tab:soilparams_mms}
\begin{tabular}{|c|c|c|c|c|c|}
\hline
$\theta_s$ & $\theta_r$ & $n$ & $K_s$ [m/day] & $\alpha$ [m$^{-1}$] & $h_{\text{cap}}$ [m] \\
\hline
0.41 & 0.095 & 2.0 & 0.0624 & 1.9 & $\alpha^{-1}$ \\
\hline
\end{tabular}
\end{table}

Instead of prescribing a physically motivated initial state, we impose a
manufactured solution $u_{\text{exact}}(x,t)$ that is smooth in both space and
time. This ensures that the exact solution is available in closed form,
allowing us to assess the accuracy of the numerical scheme and compute
convergence rates. 

The boundary conditions are imposed consistently with the manufactured
solution:
\[
u(x,t) = u_{\text{exact}}(x,t), \qquad x=0,1, \,t \in [0,1].
\]
and the initial condition is given by:
$$
u(x,0) = u_{exact}(x,0), \quad x\in [0,1].
$$
The manufactured solution is
\begin{equation}
    u_{exact}(x,t) = \frac{1}{2} m \mathcal{B}(1,m,1/n) \left[ \tanh\!\left(20\frac{x - (0.5-0.25t)}{\tanh(10 t) + 0.1}\right) + 1 \right].
\end{equation}

Here, \(\mathcal{B}\) denotes the Beta function. 

From this exact solution, we compute the corresponding forcing term. 
It is obtained by substituting \(u_{exact}\) into the Richards equation. 
The resulting source term is added to the PDE, which ensures that \(u_{exact}\) is an exact solution of the modified equation. 
This procedure is standard in the method of manufactured solutions. 
It enables rigorous convergence studies. 

\begin{figure}[H]
    \centering
    \includegraphics[width=\linewidth]{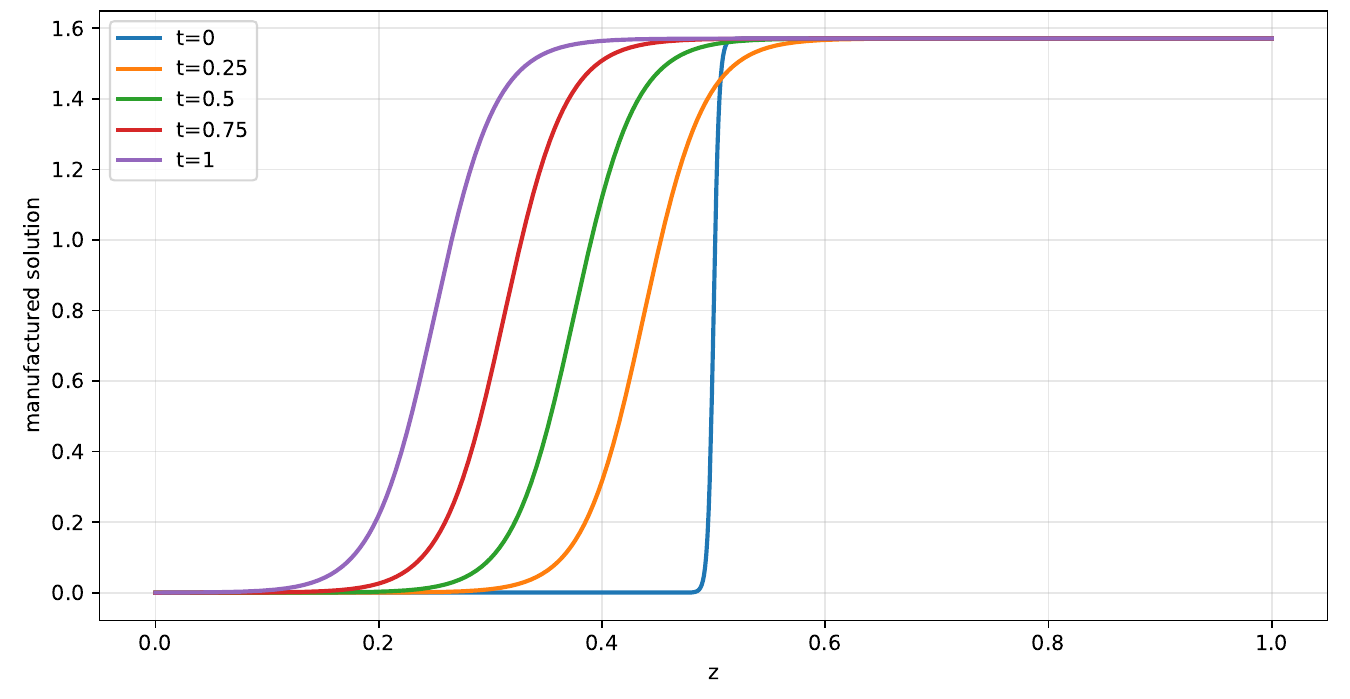}
    \caption{Manufactured solution at selected times.}
    \label{fig:manufactured}
\end{figure}

\begin{figure}[H]
    \centering
    \includegraphics[width=1\linewidth]{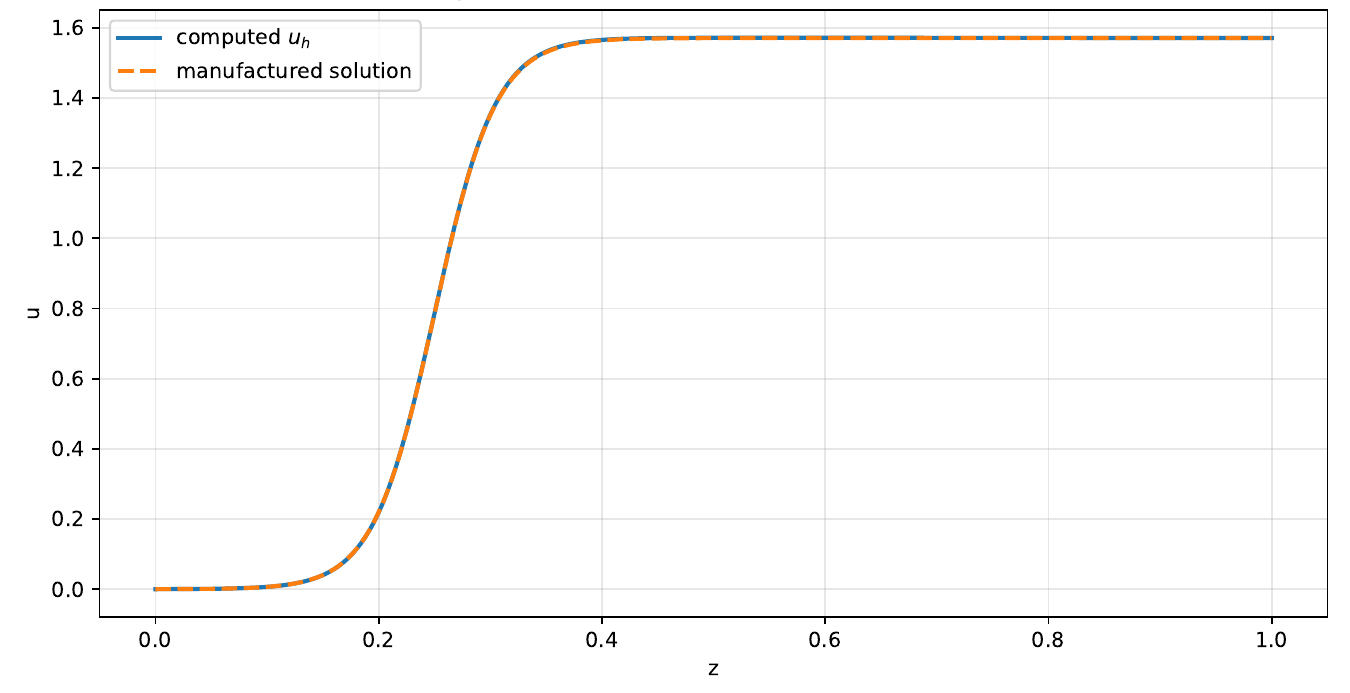}
    \caption{Computed solution vs. manufactured solution at final time ($t=1$ day).}
    \label{fig:manufactured-compare}
\end{figure}

Figure~\ref{fig:manufactured} presents the manufactured solution at several selected times. 
The plots illustrate the propagation and smoothing behaviour. 
A comparison at the final time is given in Figure~\ref{fig:manufactured-compare}. 
It shows excellent agreement between the numerical and manufactured solutions.

A spatial convergence study is carried out with a fixed time step \(\Delta t = 5 \times 10^{-5}\). 
The results, summarized in Table~\ref{tab:space_conv}, show an order in space approaching \(p \approx 2.0\) on the finest meshes. 
The reported time corresponds to the wall-clock time.


\begin{table}[H]
\centering
\caption{Space convergence (fixed $\Delta t = 5 \times 10^{-5}$).}
\label{tab:space_conv}
\begin{tabular}{c c c c}
\hline
$h$ & $L^2$ error & p & time (s) \\
\hline
$1.42857143 \times 10^{-2}$ & $6.037751 \times 10^{-3}$ & - & 25.1 \\
$7.14285714 \times 10^{-3}$ & $8.988816 \times 10^{-4}$ & $2.74781$ & 47.2 \\
$3.57142857 \times 10^{-3}$ & $1.389318 \times 10^{-4}$ & $2.69375$ & 82.8 \\
$1.78571429 \times 10^{-3}$ & $3.47545 \times 10^{-5}$ & $1.99911$ & 151.0 \\
\hline
\end{tabular}
\end{table}

Similarly, the temporal convergence is examined with a fixed mesh size \(h = 1 \times 10^{-3}\). Table~\ref{tab:time_conv} shows that the order is consistently close to 1, as expected for the chosen time discretization method.

\begin{table}[H]
\centering
\caption{Time convergence (fixed $h = 1 \times 10^{-3}$).}
\label{tab:time_conv}
\begin{tabular}{c c c c}
\hline
$\Delta t$ & $L^2$ error & $q$ & time (s) \\
\hline
$1 \times 10^{-2}$ & $8.232503 \times 10^{-3}$ & - & 6.9 \\
$5 \times 10^{-3}$ & $3.895108 \times 10^{-3}$ & $1.07967$ & 11.1 \\
$2.5 \times 10^{-3}$ & $1.918802 \times 10^{-3}$ & $1.02146$ & 19.3 \\
$1.25 \times 10^{-3}$ & $9.598079 \times 10^{-4}$ & $0.999388$ & 33.3 \\
\hline
\end{tabular}
\end{table}

\subsection{2D fully saturated test case}
In this section, we conduct a 2D numerical simulation to compare the results obtained by the Hydrus software, which solves the problem in terms of \( \Psi \), and the solution obtained with the $u$-formulation. The test case presented here is based on fully saturated regions. The van Genuchten-Mualem model is employed, with parameters provided in tables \ref{tab:Leverett function} and \ref{tab:u variable}. 

\begin{figure}[H]
    \centering
    \begin{tikzpicture}[scale=0.8]
    \draw[<-, dashed] (0,2) node[left] {$z$} -- (0,-9) ;
    \draw[->, dashed] (-2,-8) -- (10,-8) node[above] {$x$}; 
    \draw[color=blue, very thick] (2,0) -- (6,0);
    \draw[->, color=blue, thick,dashed] (5,1) -- (5,0.5);
    \draw[->, color=blue, thick, dashed] (3,1) -- (3,0.5);
    \node[color=blue] at (4,1.5) {fixed water content $\theta=0.41$ on $\Gamma_D$};
    \node[color=blue] at (4,-0.5) {$\Gamma_D$};
    \draw[color=red]  (8,0) -- (8,-8);
    \draw[->, color=red, thick] (8.25,-2) -- (9,-2);
    \draw[color=red, thick] (8.5,-2.25) -- (8.75,-1.75);
    \draw[->, color=red, thick] (8.25,-1) -- (9,-1);
    \draw[color=red, thick] (8.5,-1.25) -- (8.75,-0.75);
    \draw[->, color=red, thick] (8.25,-6) -- (9,-6);
    \draw[color=red, thick] (8.5,-6.25) -- (8.75,-5.75);
    \draw[->, color=red, thick] (8.25,-7) -- (9,-7);
    \draw[color=red, thick] (8.5,-7.25) -- (8.75,-6.75);
    \node[rotate=-90, color=red] at (9.5, -4) {No water flux a.e $q_N=0$};
    \node[color=red] at (7.7, -4) {$\Gamma_N$};
    \draw[color=green] (8,-8) -- (0,-8);
    \draw[<-, color=green, thick] (2,-9) -- (2,-8.5);
    \draw[<-, color=green, thick] (1,-9) -- (1,-8.5);
    \draw[<-, color=green, thick] (6,-9) -- (6,-8.5);
    \draw[<-, color=green, thick] (7,-9) -- (7,-8.5);
    \node[color=green] at (4,-9.5) {Free drainage a.e $\nabla \Psi \cdot \boldsymbol{n}= 0$};
    \node[color=green] at (4,-7.7) {$\Gamma_F$};
    \draw[color=red]  (0,0) -- (0,-8);
    \draw[color=red]  (0,0) -- (2,0);
    \draw[color=red]  (6,0) -- (8,0);
    \draw[->, color=red, thick] (-0.25,-2) -- (-1,-2);
    \draw[color=red, thick] (-0.5,-2.25) -- (-0.75,-1.75);
    \draw[->, color=red, thick] (-0.25,-1) -- (-1,-1);
    \draw[color=red, thick] (-0.5,-1.25) -- (-0.75,-0.75);
    \draw[->, color=red, thick] (-0.25,-6) -- (-1,-6);
    \draw[color=red, thick] (-0.5,-6.25) -- (-0.75,-5.75);
    \draw[->, color=red, thick] (-0.25,-7) -- (-1,-7);
    \draw[color=red, thick] (-0.5,-7.25) -- (-0.75,-6.75);
    \node[rotate=90, color=red] at (-1.5, -4) {No water flux a.e $q_N=0$};
    \node[color=red] at (0.5, -4) {$\Gamma_N$};
    \node at (4,-4) {$\Omega$};
    \node at (4,-3) {$\theta=0.13$ at $t=0$.};
\end{tikzpicture}
    \caption{The Richards equation boundary conditions}
    \label{fig:2dBCmiddle}
\end{figure}
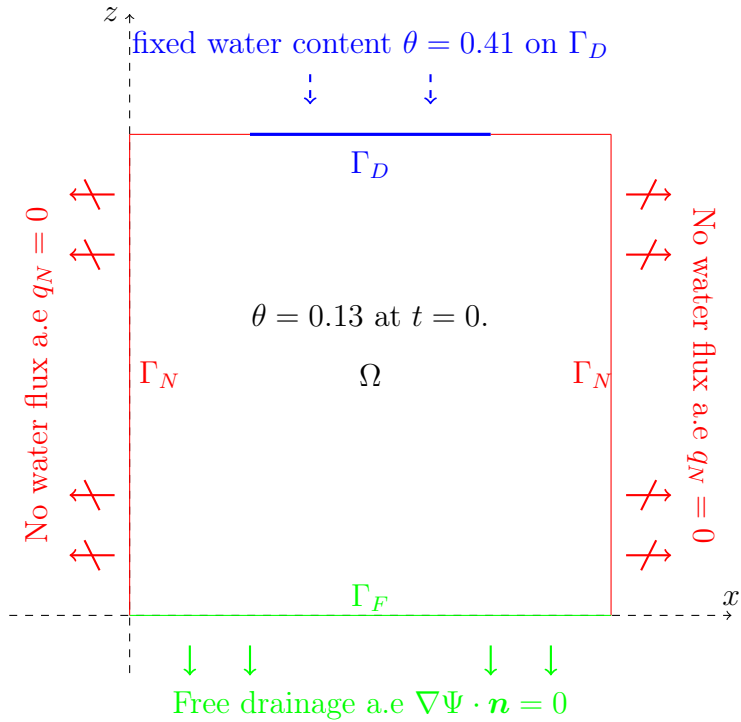

The domain is given by $\Omega=(0, 100)^2$, with spatial coordinates measured in centimeters. The top boundary consists of a Dirichlet segment \(\Gamma_D\) of length 50 cm, centered horizontally along the top edge. The rest of the top boundary is part of \(\Gamma_N\). The left and right boundaries are part of \(\Gamma_N\) and the bottom boundary is part of \(\Gamma_F\). The initial and boundary conditions are defined as follows
\begin{align*}
    &\theta(\boldsymbol{x},0) = 0.13,\;\forall \boldsymbol{x} \in \Omega, \\
    &\theta_D = 0.41,\\
    &q_N = 0,\\
    &\nabla \Psi \cdot \boldsymbol{n}=0, \text{ on } \Gamma_F.
\end{align*}
Figure \ref{fig:2dBCmiddle} gives an illustration of boundary and initial values. Soil parameters are given in table \ref{tab:2dtestparam}.

\begin{table}[H]
    \centering
    \caption{Values of soil hydraulic parameters}
    \begin{tabular}{|c|c|c|c|c|c|}
         \hline$\theta_s$ & $\theta_r$ & $n$ & $K_s\;[\mathrm{cm/h}]$ & $\alpha\;[\mathrm{cm^{-1}}]$ & $h_{cap}\;[\mathrm{cm}]$\\
\hline 0.41 & 0.047  & 1.48 & 1.96 & 0.015 & $ \alpha^{-1}$\\
\hline
    \end{tabular}
    
    \label{tab:2dtestparam}
\end{table}

\begin{figure}[H]
    \centering
    \includegraphics[width=1\linewidth]{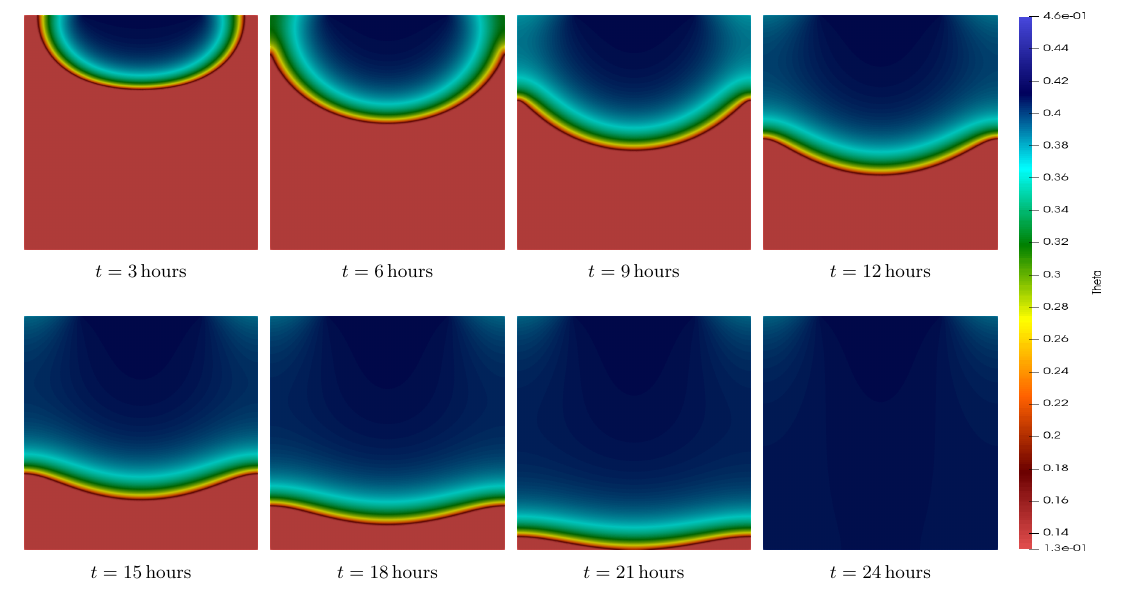}
    \caption{Color plot of $\theta$ over different times using $u$-formulation.}
    \label{fig:middle2dtest}
\end{figure}

\begin{figure}[H]
    \centering
    \includegraphics[width=1\linewidth]{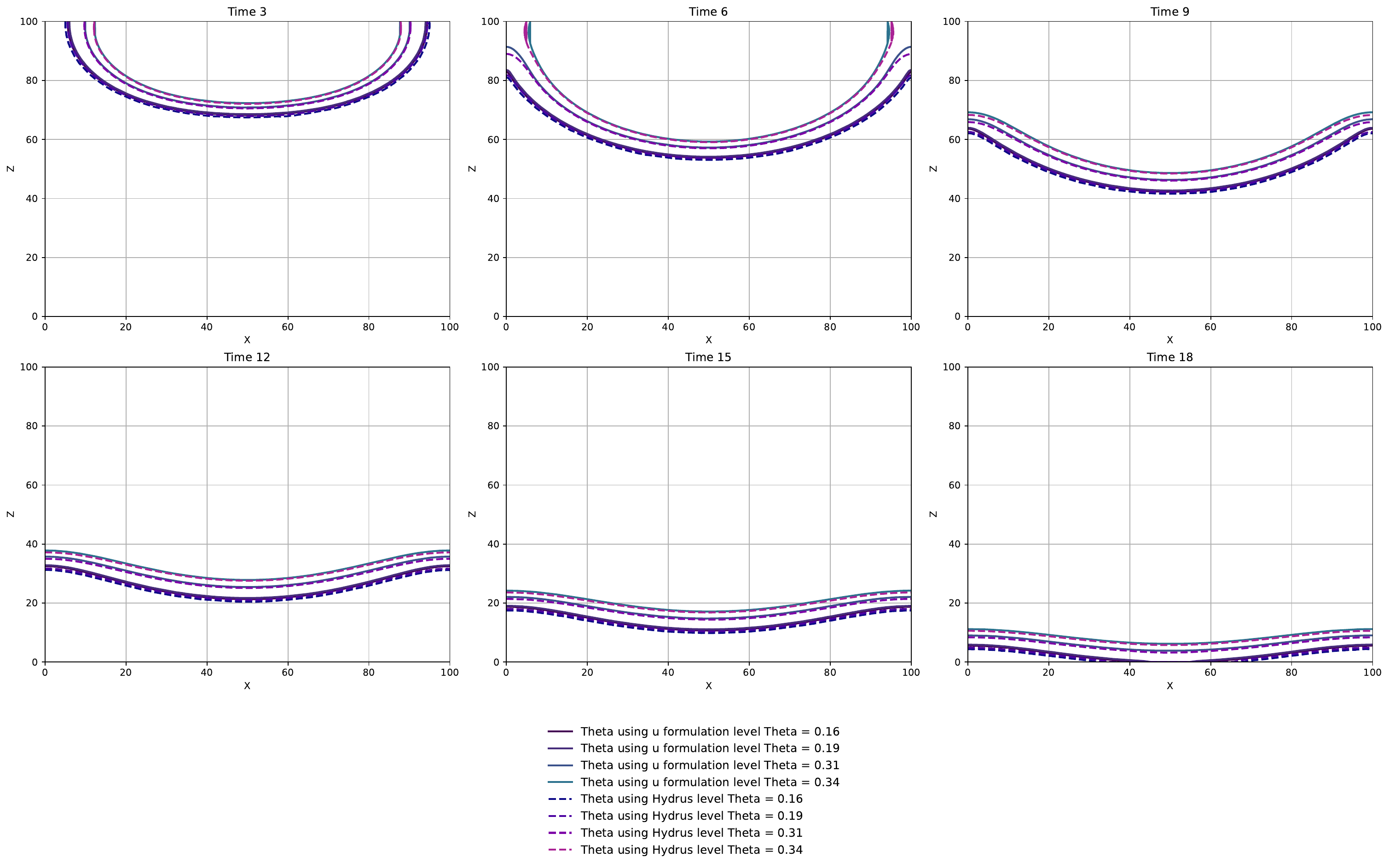}
    \caption{Contours for saturated flow.}
    \label{fig:solvshydrus2d}
\end{figure}

Figure \ref{fig:middle2dtest} shows the evolution of \( \theta \) at different times using the \( u \)-formulation. The simulation was performed on a mesh with 200 vertices per boundary segment, resulting in 95,142 triangular elements and 47,972 vertices. The time step used is \( \Delta t = 0.001 \;\mathrm{hours} \). Here, we noticed numerically that the linear semi-implicit method is sufficient to solve the $u$-formulation. The color plots represent snapshots of \( \theta \) at various time intervals, showing the progression of the solution over time. Figure \ref{fig:solvshydrus2d} illustrates the evolution of \( \theta \) over different times, where we compare the contours of the solution obtained by the Hydrus software and the proposed \( u \) formulation using the same numerical parameters as in Figure \ref{fig:middle2dtest}. Our results closely match those from the Hydrus software, confirming that the proposed numerical method produces accurate and reliable solutions.

\subsection{Layered soil test case}
\label{subsec:layredsoil}
\label{Hetsatunsattestcase}
\noindent In this section, we solve the 2D layered-soil benchmark used in \cite{Manzini2004} using a domain-decomposition approach combined with our new formulation. The square domain 
\[
\Omega = (0 ,100)^2
\]
with spatial coordinates measured in centimeters, is partitioned along the curved interface \(\Gamma\) (see figure \ref{fig:layeredsoil}), which separates two soil layers with distinct hydraulic properties. The hydraulic functions are modeled by the van Genuchten equations for the capillary pressure and the van Genuchten–Mualem model for the relative permeability. We define the subdomains by the curve 
\[
\Gamma = \{(x,\xi(x)):\; x\in [0,100]\},
\]
and we define $\xi$ by
\[
\xi(x) = 100 \!\left( 0.1 \bigl(1 - \cos(\tfrac{\pi x}{100})\bigr) + 0.45 \right).
\]

The soil parameters are given in table \ref{tab:parameters}.

\begin{table}[H]
    \centering
    \caption{Soil parameters for different regions in the layered soil.}
    \renewcommand{\arraystretch}{1.2}
    \begin{tabular}{|c |c |c |c |c |c|}
        \hline
        Index & $\theta_s$ & $\theta_r$ & $\alpha$ (cm$^{-1}$) & $n$ & $K_s$ (cm$\cdot$h$^{-1}$) \\
        \hline
        1 (Top)& 0.50 & 0.120 & 0.028 & 3.00 & 0.25 \\
        \hline
        2 (Bottom)& 0.46 & 0.034 & 0.016 & 1.37 & 2.00 \\
        \hline
    \end{tabular}
    \label{tab:parameters}
\end{table}
The initial condition is given by $\Psi(x,z,0)=-z$. Fully saturated Dirichlet boundary conditions are applied on the top and bottom parts, i.e. $S=1$, and no flux Neumann conditions on the vertical parts of the domain boundary.

\begin{figure}[H]
    \centering
    \includegraphics[width=1\linewidth]{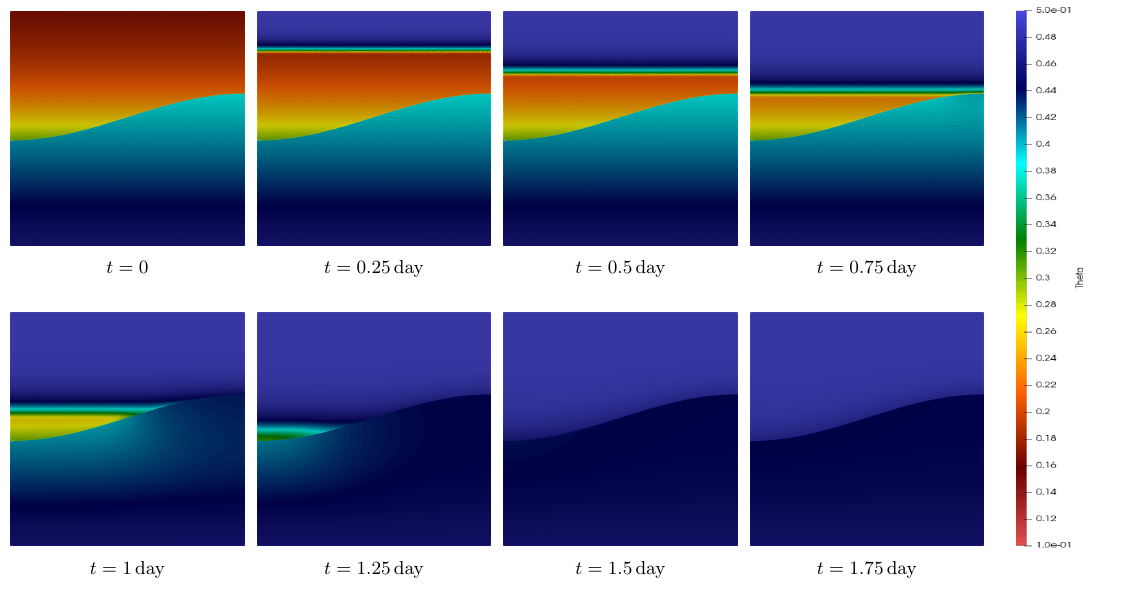}
    \caption{Water content $\theta$ over various times for Manzini test case.}
    \label{fig:hetres}
\end{figure}

\begin{figure}[H]
    \centering
    \includegraphics[width=1\linewidth]{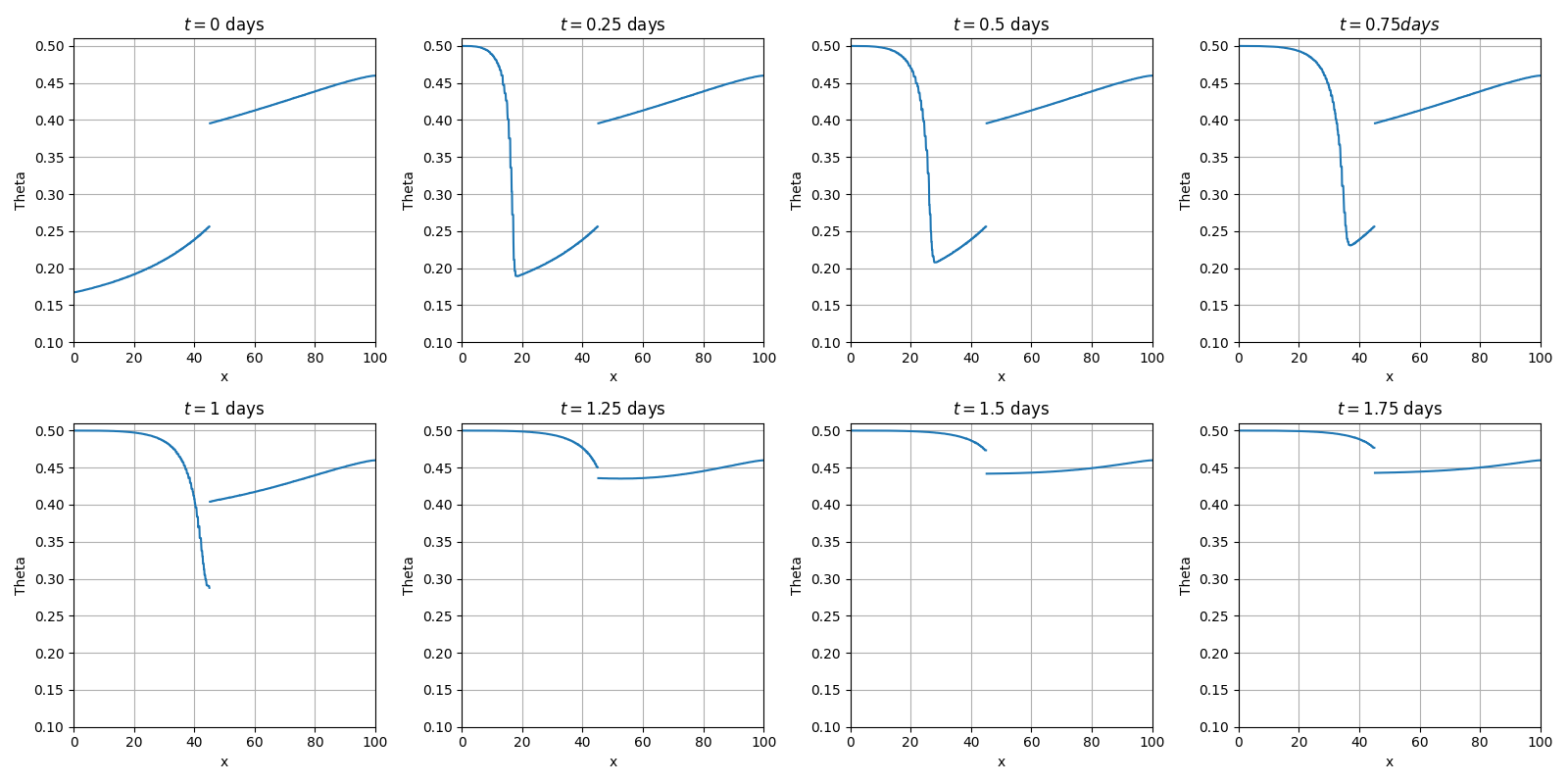}
    \caption{Plot over a vertical cut at $x=50$ from figure \ref{fig:hetres}.}
    \label{fig:hetres1d}
\end{figure}
We solve this problem using a mesh with 200 vertices on each segment of the boundary and on \(\Gamma\), resulting in 89,306 triangular elements and 45,064 vertices. The time step is set to \(\Delta t = 0.05\,\) hours. A tolerance of \(10^{-3}\) is used with a maximum of 10 iterations, with $\lambda=25$. The results shown in Figures \ref{fig:hetres} and \ref{fig:hetres1d} are in good agreement with those reported in~\cite{Manzini2004}. 

\subsection{Fully saturated and unsaturated layered soil test case}
To further validate the proposed method, in this section we present an additional numerical test in the layered soil test case of Section~\ref{subsec:layredsoil}, with a change to the initial condition, which is given by:
\begin{equation*}
    S(x,z,0)=
    \begin{cases}
        1, & \text{if } z \leq \xi(x),\\
        0, & \text{if } z > \xi(x).
    \end{cases}
\end{equation*}
This test case demonstrates that our method can correctly resolve such a sharp interface, with fully saturated soil on one side and fully dry soil on the other side.

\begin{figure}[H]
    \centering
    \includegraphics[width=1\linewidth]{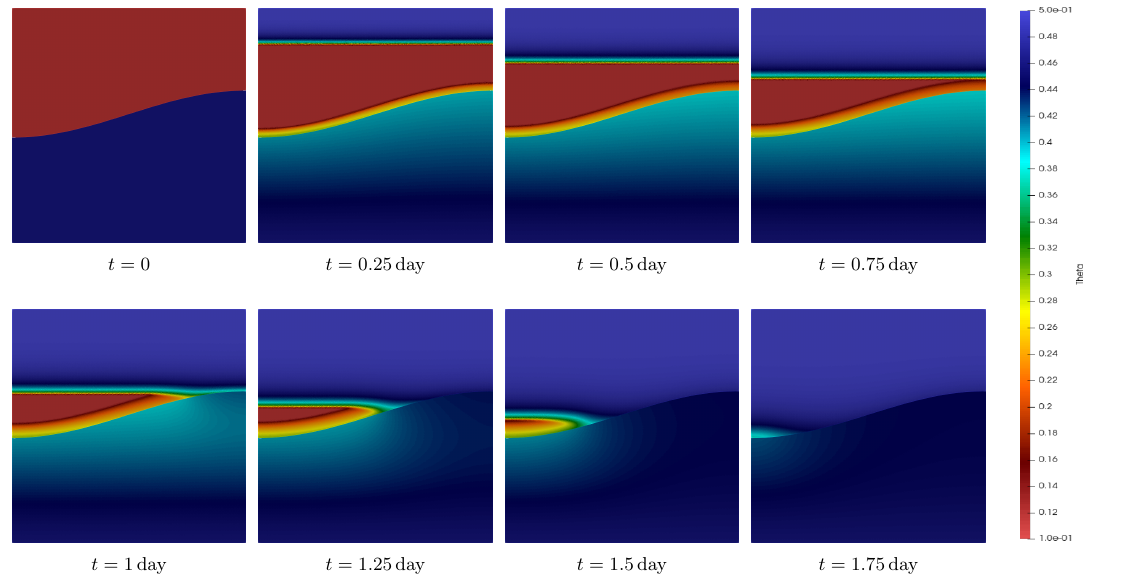}
    \caption{Water content $\theta$ at various times.}
    \label{fig:heterpic}
\end{figure}

\begin{figure}[H]
    \centering
    \includegraphics[width=1\linewidth]{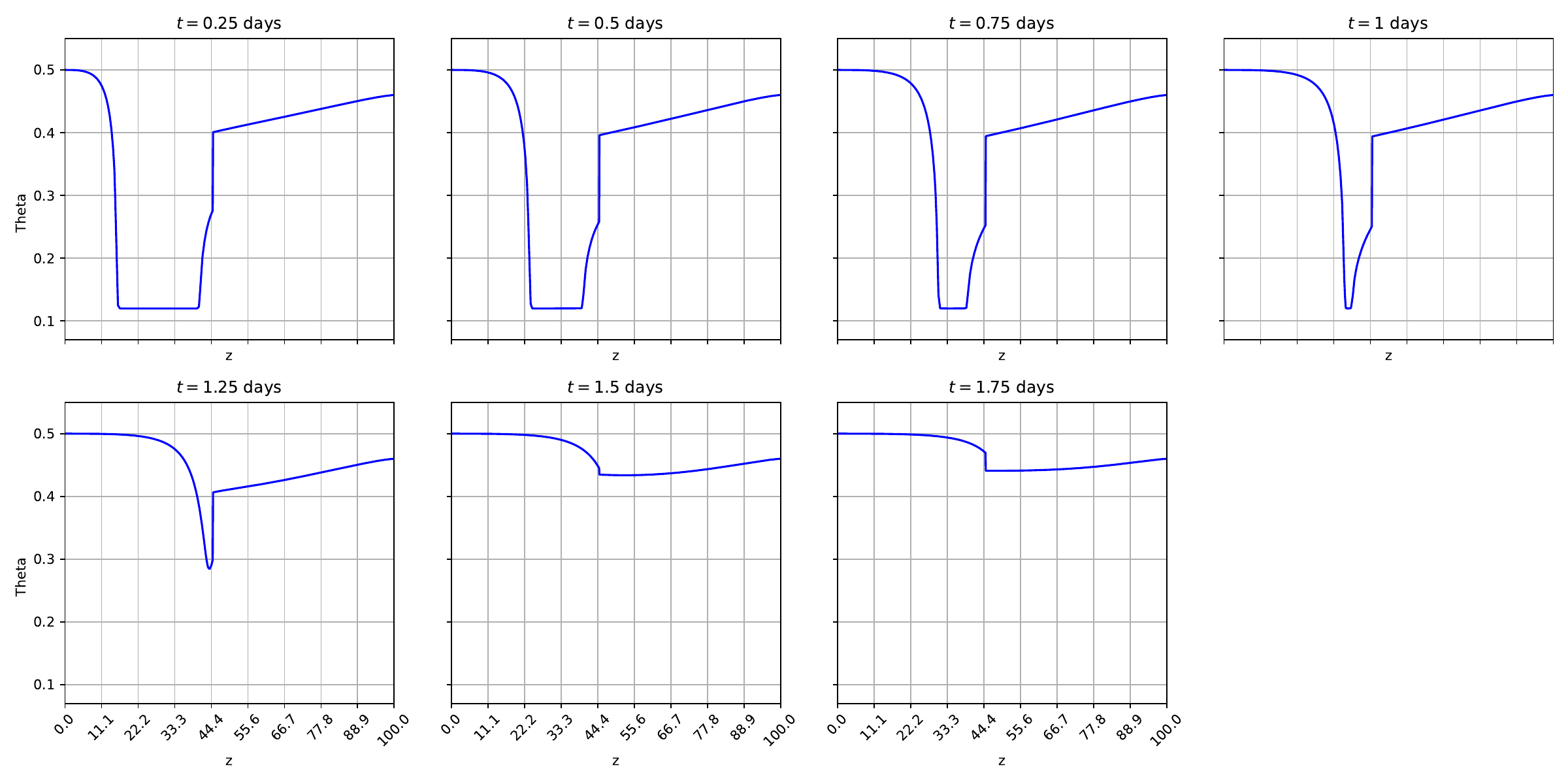}
    \caption{Plot over a vertical cut at $x=50$ from figure \ref{fig:heterpic}.}
    \label{fig:heterpic1d}
\end{figure}

Figures~\ref{fig:heterpic} and \ref{fig:heterpic1d} show the evolution of water content \(\theta\) over time. The discontinuity across the horizontal interface reflects the contrast in soil properties: water infiltrates into the upper layer due to the capillary forces visible in Figure~\ref{fig:heterpic}.

To show that $\Psi$ remains continuous even in almost dry regions, we repeat the test with a modified initial value to avoid the degeneracy of the pressure head (i.e., $\Psi \to \infty$ as $S\to 0$). The initial condition is
\[
S(x,z,0)=
\begin{cases}
1, & z\le \xi(x),\\
\epsilon, & z>\xi(x).
\end{cases}
\]

\begin{figure}[H]
    \centering
    \includegraphics[width=1\linewidth]{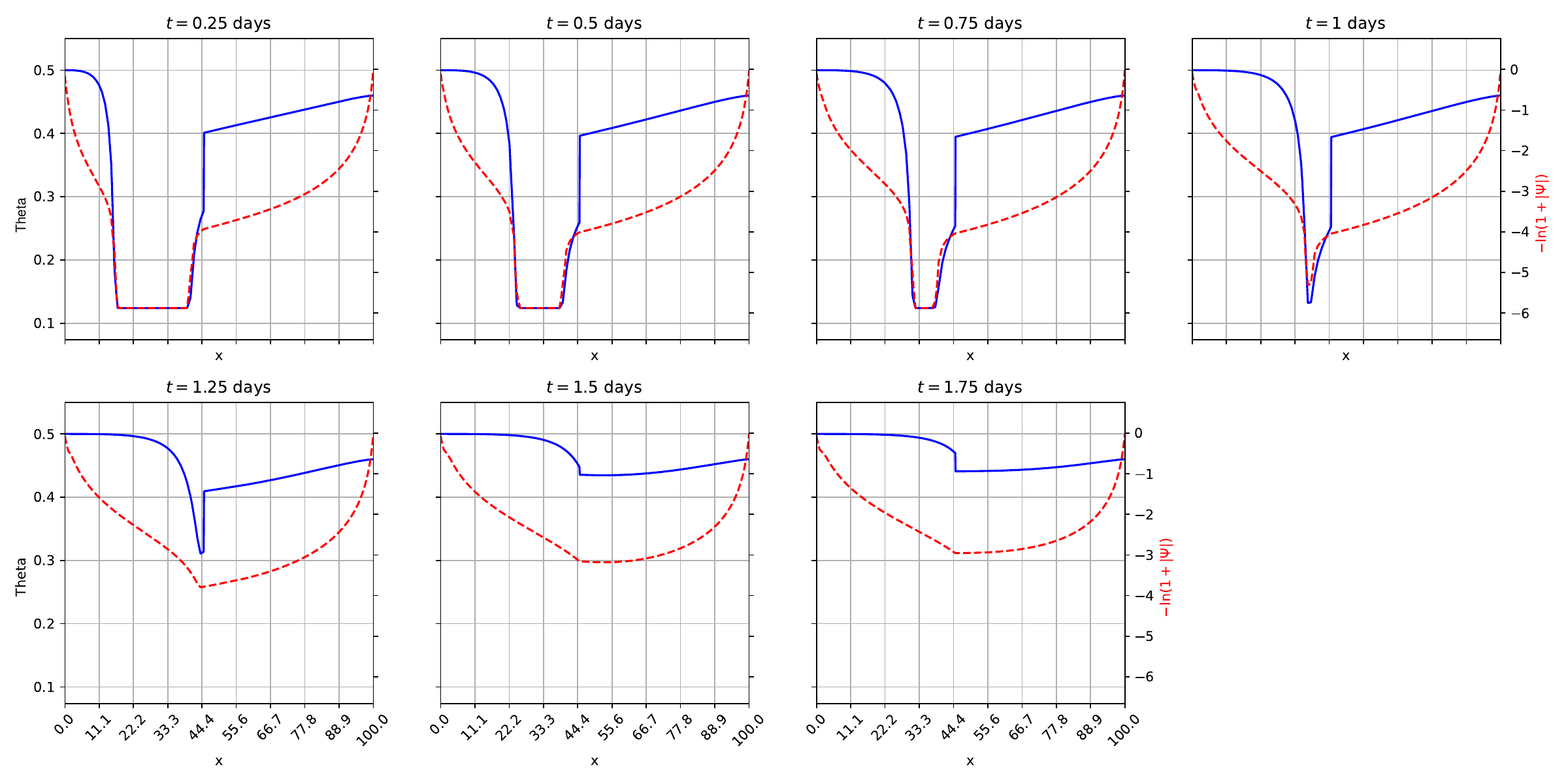}
    \caption{Plot over a line at $x=50$ for $\epsilon = 10^{-2}$ comparing $\theta$ and $-\ln(1+|\Psi|)$.}
    \label{fig:plotoverlineeps1e2}
\end{figure}

\begin{figure}[H]
    \centering
    \includegraphics[width=1\linewidth]{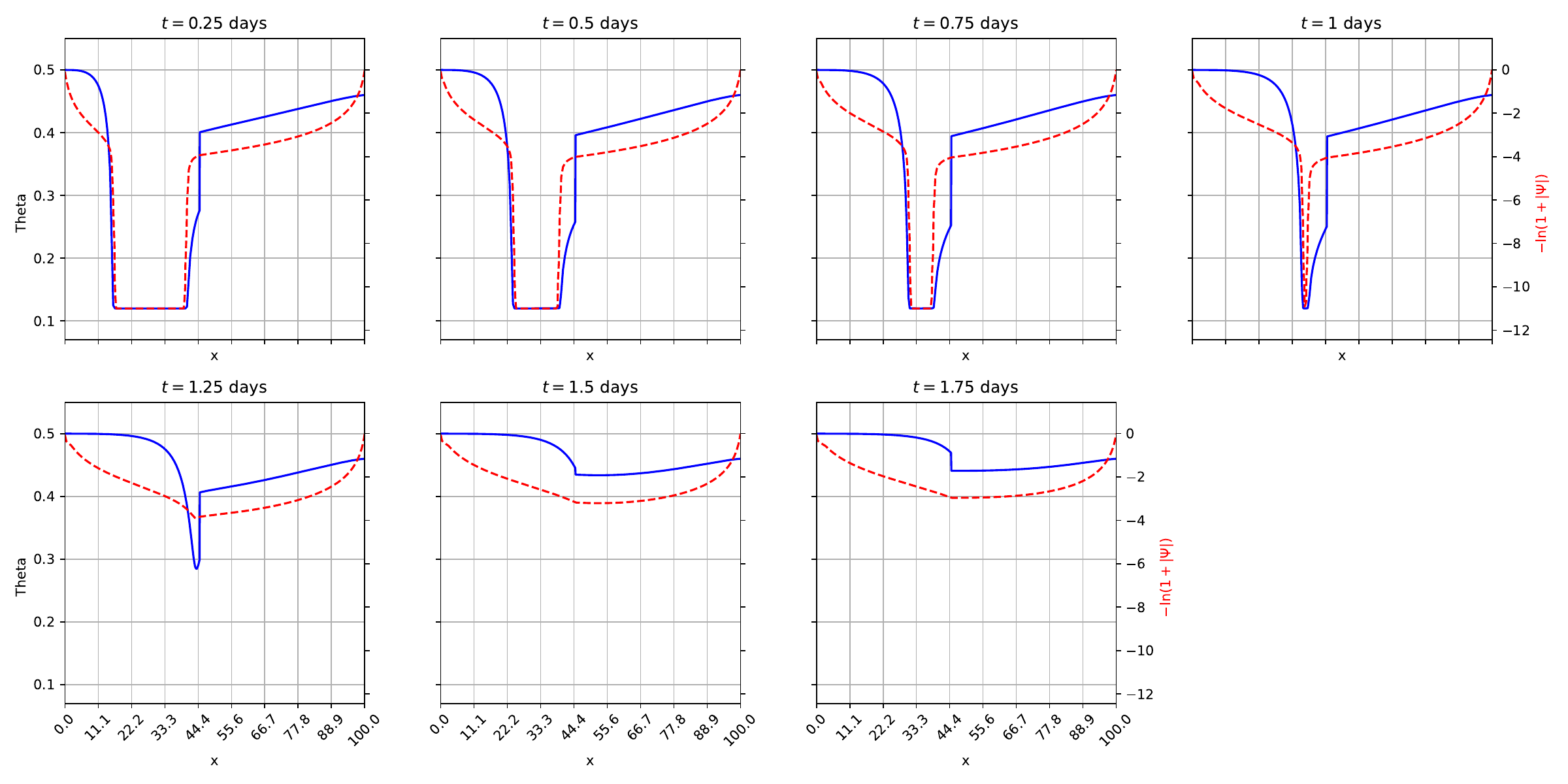}
    \caption{Plot over a line at $x=50$ for $\epsilon = 10^{-8}$ comparing $\theta$ and $-\ln(1+|\Psi|)$.}
    \label{fig:plotoverlineeps1e8}
\end{figure}

Here, $\epsilon>0$ is a small value. We have solved the equations using the same parameters as in section~\ref{subsec:layredsoil}. Figures \ref{fig:plotoverlineeps1e2} and \ref{fig:plotoverlineeps1e8} shows the values of \( \theta \) (in blue) and \( -\ln(1 + |\Psi|) \) (in red) along a vertical line at \( x = 50 \), for two different initial values: \( \epsilon = 10^{-2} \) (Figure \ref{fig:plotoverlineeps1e2}) and \( \epsilon = 10^{-8} \) (Figure \ref{fig:plotoverlineeps1e8}). Despite the sharp transitions in \( \theta \), the plots of \( -\ln(1 + |\Psi|) \) remain smooth, confirming that \( \Psi \) stays continuous even for very small values of \( \epsilon \), which shows that the non-overlapping domain decomposition method used performs well.

\subsection{Heterogeneous test case: a soil inclusion}

To further illustrate the flexibility of the domain decomposition approach, we consider a setting where one soil region is fully embedded within another, rather than the layered configuration of the previous tests. The hydraulic parameters are kept the same as in Section~\ref{Hetsatunsattestcase}. The computational domain is defined as
\[
\Omega_1 = [0,100]^2 \setminus \overline{\Omega}_2, \qquad 
\Omega_2 = \left\{ (x,y)\in [0,100]^2 : (x-50)^2 + (y-50)^2 < 25^2\right\}.
\]

No-flux boundary conditions are imposed on $\partial \Omega$. 
The initial conditions are set as $S = \epsilon$ in $\Omega_1$ and $S = 1$ in $\Omega_2$, with $\epsilon \geq 0$. 
The problem is discretized using 21{,}956 triangular elements and 11{,}179 vertices. 
We employ a time step of $\Delta t = 0.1$ hours and choose $\lambda = 25$. 
On average, the nonlinear solver requires 4 iterations to converge, with an absolute tolerance of $10^{-3}$.  

The $\Psi$-formulation is solved using the modified L-scheme \cite{Mitra2019}, 
with $\mathcal{M}\Delta t = 2 \times 10^{-5}$. 
The hydraulic conductivity term $K_r$ is lagged in the same way as in the semi-implicit method \eqref{eq:semiimplicitscheme}. 
Calculations are performed on the same mesh and with the same time step, 
and the nonlinear solver converges in about 8 iterations on average, 
with an absolute tolerance of $10^{-3}$.

\begin{figure}[H]
    \centering
    \includegraphics[width=1\linewidth]{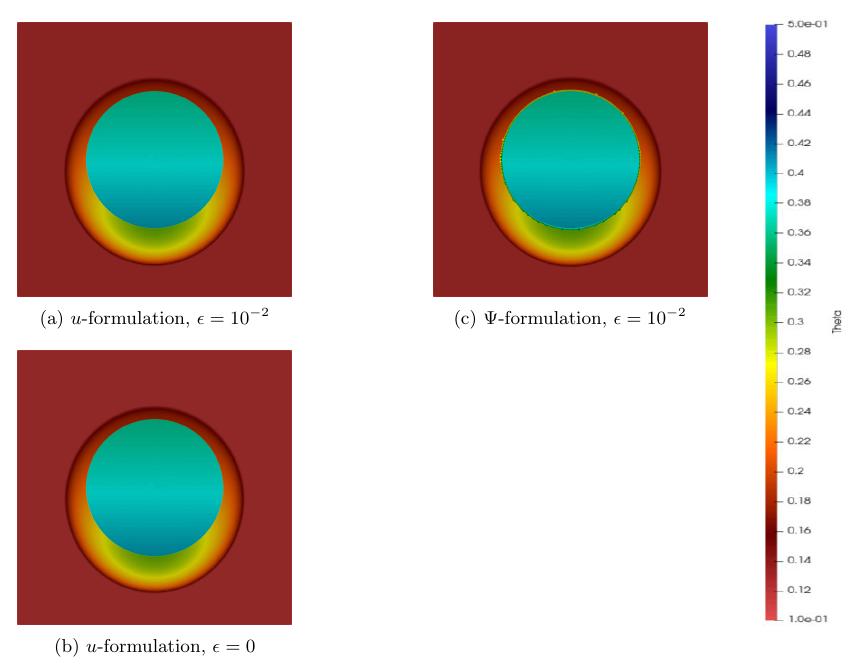}
    \caption{Numerical comparison of water content $\theta$: $u$-formulation vs. $\Psi$-formulation at final time $t=42$ hours.}
    \label{fig:hetinside}
\end{figure}

As shown in Figure~\ref{fig:hetinside}, both the $u$-formulation and the $\Psi$-formulation provide
similar results for the water content $\theta$ at the final simulation time for $\epsilon=10^{-2}$, and the figure also shows the results obtained by the $u$-formulation for a completely dry medium ($\epsilon=0$). 

For smaller values of $\epsilon$, the convergence of the nonlinear solver for the $\Psi$-formulation deteriorates significantly. For example, with $\epsilon = 10^{-5}$, the solver requires more than $1000$ iterations in the first time step, yet the residual tolerance remains at approximately $0.04$.
In contrast, the $u$-formulation remains stable for $\epsilon=0$ and converges in about 4 iterations.

\begin{figure}[H]
    \centering
    \includegraphics[width=\linewidth]{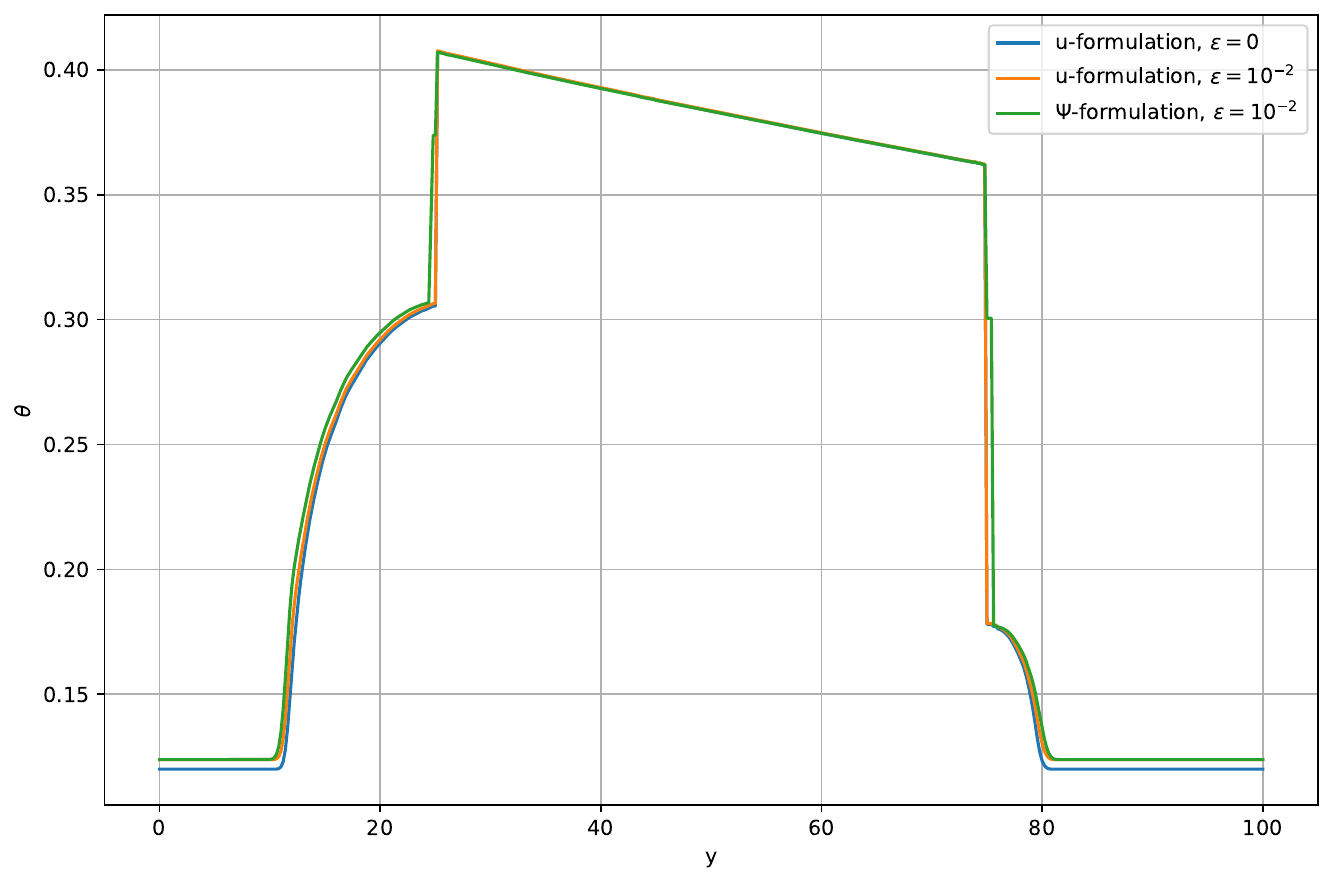}
    \caption{Water content $\theta$ on the vertical cut $x=50$ cm and final simulation time, 
    comparing the $u$-formulation (for $\epsilon=0$ and $\epsilon=10^{-2}$) with the $\Psi$-formulation (for $\epsilon=10^{-2}$).}
    \label{fig:theta_line}
\end{figure}

A more detailed comparison of the water content distribution along the vertical profile at 
$x=50$ cm is shown in Figure~\ref{fig:theta_line}. The line plots confirm the close agreement 
between the two formulations for $\epsilon=10^{-2}$, while also showing the line plot for the result obtained by the $u$-formulation for $\epsilon=0$.

\section{Conclusion}\label{sec:conclusion}
We developed a finite element method for the Richards equation based on a new bounded auxiliary variable $u$, which removes unbounded terms from the weak formulation. The method uses a semi-implicit discretization with Newton iterations and we extend it to heterogeneous media through a non-overlapping Schwarz domain decomposition. Numerical tests with the van Genuchten model, including dry fibrous sheets, mixed saturated/unsaturated regions, manufactured solutions, and layered or heterogeneous soils, demonstrate stability and accuracy across the full range $S \in [0,1]$ without variable switching or regularization. These results establish the $u$-formulation as a robust alternative to classical approaches. Future work will focus on rigorous numerical analysis and on coupling with solute and heat transport to address broader multiphysics applications.

\section*{Acknowledgments}

This work was supported by an NSERC, Canada Discovery Grant (RGPIN-2019-06855) to Yves Bourgault 
and an NSERC, Canada Discovery Grant (RGPIN/5220-2022 \& DGECR/526-2022) to Abdelaziz Beljadid.



\biboptions{sort&compress} 
\bibliographystyle{elsarticle-num}
\bibliography{references}

\end{document}